\documentclass[twocolumn]{autart} 

\usepackage{graphicx}
\usepackage{amsmath}
\usepackage{amssymb}
\usepackage{mathtools} 
\usepackage[utf8]{inputenc}
\usepackage{letltxmacro}
\usepackage[table,xcdraw]{xcolor}
\usepackage{hhline} 

\newcommand{\nref}[1]{(\ref{#1})}
\newcommand{\bfs}[1]{\boldsymbol{#1}}
\newcommand{\col}[1]{\operatorname{col}\left({#1}\right)}
\newcommand{\diag}[1]{\operatorname{diag}\left({#1}\right)}
\newcommand{\kraj}{\hspace*{\fill} \qed}
\newcommand{\krajdokaz}{\hfill $\blacksquare$}

\newcommand{\z}[1]{\left( #1 \right)}
\newcommand{\zs}[1]{\left[ #1 \right]}

\newcommand{\sm}[1]{\left[\begin{smallmatrix} #1 \end{smallmatrix}\right]}
\newcommand{\n}[1]{\left\| #1 \right\|}
\newcommand{\red}{\nonumber \\}
\newcommand{\vprod}[2]{\left \langle #1 \ \middle\vert\  #2 \right \rangle}
\newcommand{\oln}[1]{\overline{#1}}

\newcommand{\R}{\mathbb{R}}
\newcommand{\Z}{\mathbb{Z}}
\newcommand{\N}{\mathbb{N}}

\newcommand{\dom}{\mathrm{dom}}
\newcommand{\proj}{\mathrm{proj}}

\usepackage[utf8]{inputenc}
\usepackage[shortlabels]{enumitem} 
\usepackage{cuted}
\usepackage{widetext}

\newcommand{\Diag}[1]{\operatorname{blkdiag}\left({#1}\right)}

\newtheorem{assum}{Assumption}
\newtheorem{theorem}{Theorem}
\newtheorem{lemma}{Lemma}

\newtheorem{corollary}{Corollary}
\newtheorem{defn}{Definition}
\allowdisplaybreaks 
\begin{document}

\begin{frontmatter}
\runtitle{Insert a suggested running title}  

\title{A discrete-time averaging theorem and its application to zeroth-order Nash equilibrium seeking} 

\thanks[footnoteinfo]{This work was partially supported by the ERC under research project COSMOS (802348). E-mail addresses: \texttt{\{s.krilasevic-1, s.grammatico\}@tudelft.nl}.}

\author{Suad Krila\v sevi\' c} and
\author{\ Sergio Grammatico} 
          
\address{Delft Center for Systems and Control, TU Delft, The Netherlands}


\begin{abstract}                          
In this paper we present an averaging technique applicable to the design of zeroth-order Nash equilibrium seeking algorithms. First, we propose a multi-timescale discrete-time averaging theorem that requires only that the equilibrium is  semi-globally practically stabilized by the averaged system, while also allowing the averaged system to depend on ``fast" states. Furthermore, sequential application of the theorem is possible, which enables its use for multi-layer algorithm design. Second, we apply the aforementioned averaging theorem to prove semi-global practical convergence of the zeroth-order information variant of the discrete-time projected pseudogradient descent algorithm, in the context of strongly monotone, constrained Nash equilibrium problems. Third, we use the averaging theory to prove the semi-global practical convergence of the asynchronous pseudogradient descent algorithm to solve strongly monotone unconstrained Nash equilibrium problems. Lastly, we apply the proposed asynchronous algorithm to the connectivity control problem in multi-agent systems.
\end{abstract}

\end{frontmatter}

\section{Introduction}

Given a complex dynamical system, averaging techniques are used to construct a simpler system, called the \textit{averaged  system}, that is easier to analyze than the given one. Ideally, the averaged system should satisfy certain properties so that it is possible to infer stability properties of the original system based on the averaged one. These techniques are used extensively in extremum seeking results, in continuous-time systems \cite{krstic2000stability}, \cite{frihauf2011nash},  \cite{poveda2020fixed}, discrete-time systems \cite{choi2002extremum}, \cite{stankovic2011distributed}, and hybrid systems \cite{poveda2017framework}, \cite{poveda2021robust}.  \\ \\

\emph{Literature review:} Discrete-time averaging techniques have received intense attention over the years. In \cite{bai1988averaging}, the authors show that the original dynamics render the equilibrium exponentially stable under the assumption of exponential stability of the equilibrium for the averaged dynamics. Furthermore, they prove a similar result with similar assumptions for a mixed time-scale system where the fast dynamics converge to zero. The requirement for exponential stability is relaxed in \cite{wang2009input} to just semi-global practical asymptotic stability, for single time-scale systems. Furthermore, the authors include noise into their analysis and provide input-to-state stability results. In \cite{yang2022periodic}, the authors provide upper bounds for the time-scale separation parameter in the case of linear switched averaged systems by using a time-delay approach and similar assumptions as in \cite{bai1988averaging}. The previously mentioned single time-scale results assume that the jump mapping is time-varying and that this dependence gets ``smoothed out" using the averaging technique. Thus, the main source of ``perturbations" in the original system is this time dependence. Likewise, it is possible to assume that jump mapping is a function of some stochastic perturbations and that the goal of the averaging is to ``smooth out" the dependence on these perturbations. In \cite{liu2015stochastic}, the authors prove that under certain technical assumptions, the discrete-time stochastic algorithm can be approximated by its continuous counterpart, and that equilibrium of the original dynamics is weakly exponentially stable if the equilibrium of the continuous counterpart is exponentially stable. \\
The usual approach to design of extremum seeking algorithms consists of choosing a well-behaved full-information gradient-based algorithm in the case of optimization, or pseudogradient-based in the case of games, integrated with a (pseudo)gradient zeroth-order information estimation scheme \cite{krstic2000stability}, \cite{frihauf2011nash}, \cite{poveda2021robust}, \cite{poveda2021nonsmooth}. The produced estimate then replaces the real value of the (pseudo)gradient in the algorithm. A typical estimation technique it that of injecting sinusoidal perturbations into the inputs of a cost function, whose output is then correlated with the same perturbations. Via averaging techniques, it can be proven that this estimation behaves as the (pseudo)gradient, on average. The theory of averaging and singular perturbations for continuous and hybrid systems \cite{sanfelice2011singular}, \cite{wang2012analysis} enables the adaptation of a wide spectrum of algorithms. In \cite{krstic2000stability}, the authors adapt the gradient descent algorithm for the zeroth-order information case, together with the additional high-pass and low-pass filters to improve performance. An extremum seeking variant of the pseudogradient descent algorithm used for solving unconstrained games is presented in \cite{frihauf2011nash}. Recently, the authors in  \cite{poveda2021nonsmooth} propose a fixed-time zeroth-order algorithm for solving games, based on a similar full-information fixed-time algorithm. An accelerated first-order algorithm has been adapted for optimization problems in \cite{poveda2021robust}. Unfortunately, the same variety of extremum seeking algorithms in not available in discrete-time due to the limitations of the discrete-time averaging theory. In \cite{choi2002extremum}, the authors prove exponential convergence to the optimum of a quadratic function under the zeroth-order variant of the gradient descent algorithm with filtering. The authors in \cite{zargarzadeh2014extremum} prove ultimate boundness in a similar setup where the plant is assumed to be general dynamic nonlinear and the trajectories of the averaged system ultimately bounded. A similar approach is used in \cite{stankovic2011distributed}  to prove convergence to the Nash equilibrium in a game without constraints. In \cite{liu2015stochastic}, the authors prove stability of its stochastic variant. \\
On the other hand, zeroth-order methods that use other approaches for gradient estimation appear to be more successful and a recent overview for methods in optimization can be found here \cite{liu2020primer}. The authors in \cite{pang2022nash} solve an N-coalition game without local constraints for strongly monotone games by using Gaussian smoothing to estimate the pseudogradient, while the authors in \cite{tang2023zeroth} propose an algorithm for solving cooperative multi-agent cost minimization problem with local constraints, also with Gaussian smoothing. Both papers assume synchronous sampling of the agents, albeit with possible information delay. Similar approach to Gaussian smoothing is the residual feedback estimator that uses a previous evaluation of the cost function for the second point of the pseudogradient approximation, thus reducing the numbers of cost functions samples that need to be taken in one iteration. Using this approach, the authors in \cite{huango2023optimization} adapt two extra-gradient algorithms and prove convergence to the Nash equilibrium in pseudo-monotone plus games for diminishing step sizes and query radiuses. Authors in \cite{tatarenko2018learning} and \cite{tatarenko2020bandit} estimate the pseudogradient using the idea of continuous action-set learning automaton and prove convergence for strictly monotone games and merely monotone games, respectively, via diminishing step sizes, and Tikhonov regularization.\\
Asynchronous zeroth-order optimization algorithms have been well studied and an overview can be found here \cite{lian2016comprehensive}. For example, the authors in \cite{shen2021asynchronous} use the residual feedback estimator in an asynchronous gradient decent scheme to prove convergence. In the current state of the art, zeroth-order discrete-time Nash equilibrium seeking algorithms based on averaging use pseudogradient descent without projections, while algorithms based on other methods are more general, yet still assume synchronous sampling.\\
\emph{Contribution}: Motivated  by  the  above  literature  and open  research  problems,  to  the  best  of  our  knowledge, we consider an averaging technique for mixed time-scale discrete-time systems and merely semi-globally practically convergent averaged systems, with the application to the problem of learning Nash equilibria via zeroth-order discrete-time algorithms, in the cases of locally constrained agents, and asynchronous sampling. Specifically, our main technical contributions are summarized next:
\begin{itemize}
    \item  We extend the current results on averaging theory by using a mixed time-scale formulation of the original system and requiring that the averaged systems renders the equilibrium set SGPAS, unlike \cite[Thm. 2]{wang2009input}, where a single time-scale, time-variant system is used, and differently from  \cite[Thm. 2.2.4]{bai1988averaging}, \cite[Thm. 8.2.28]{mareels1996adaptive} where exponential stability is needed and the fast subsystem state is assumed to converge to the origin. Furthermore, we allow certain types of additive perturbation dynamics to interfere with the nominal averaging dynamics, and that the averaged jump mapping is a function of the fast states, thus enabling easier \emph{consecutive} application of the averaging theorem and the design of more complex algorithms.
    \item Enabled by our extended averaging theory, we propose two novel zeroth-order algorithms for game equilibrium seeking in discrete time. The first algorithm solves the equilibrium in games with local constraints, differently from \cite{stankovic2011distributed}, \cite{pang2022nash} where agents have no constraints; while the second one solves the problem in the case where the agents are asynchronous, i.e. the agents do not sample at the same time, nor do they coordinate in any way, differently from \cite{stankovic2011distributed}, \cite{pang2022nash} where the agents sample synchronously.
\end{itemize}

\emph{Notation}: {The set of real numbers and the set of nonnegative real numbers are denoted by $\mathbb{R}$ and  $\mathbb{R}_+$, respectively. Given a set $\mathcal{Z}$, $\mathcal{Z}^n$ denotes the Cartesian product of $n$ sets $\mathcal{Z}$.} For a matrix $A \in \mathbb{R}^{n \times m}$, $A^\top$ denotes its transpose. For vectors $x, y \in \mathbb{R}^{n}$ and $M \in \R^{n \times n}$ a positive semi-definite matrix and $\mathcal{A} \subset \R^n$, $\vprod{x}{y}$, $\|x \|$,  $\|x \|_M$ and $\|x \|_\mathcal{A}$ denote the Euclidean inner product, norm, weighted norm and distance to set respectively. Given $N$ vectors $x_1, \dots, x_N$, possibly of different dimensions, $\col{x_1, \dots x_N} \coloneqq \left[ x_1^\top, \dots, x_N^\top \right]^\top $. Collective vectors are denoted in bold, i.e,  $\bfs{x} \coloneqq \col{x_1, \dots, x_N}$ and for each $i = 1, \dots, N$, $\bfs{x}_{-i} \coloneqq \col{ x_1, \dots,  x_{i -1},  x_{i + 1}, \dots, x_N }$ { as they collect vectors from multiple agents.} Given $N$ matrices $A_1$, $A_2$, \dots, $A_N$, $\Diag{A_{1}, \ldots, A_{N}}$ denotes the block diagonal matrix with $A_i$ on its diagonal. {Given a vector $x$, $\operatorname{diag}(x)$ represents a diagonal matrix whose diagonal elements are equal to the elements of the vector $x$.} For a function $v: \mathbb{R}^{n} \times \mathbb{R}^{m}  \rightarrow \mathbb{R}$ differentiable in the first argument, we denote the partial gradient vector as $\nabla_x v(x, y) \coloneqq \col{\frac{\partial v(x, y)}{\partial x_{1}}, \ldots, \frac{\partial v(x, y)}{\partial x_{N}}} \in \mathbb{R}^{n}$. We use $\mathbb{S}^{1}:=\left\{z \in \mathbb{R}^{2}: z_{1}^{2}+z_{2}^{2}=1\right\}$ to denote the unit circle in $\R^2$. {The set-valued mapping $\text{N}_{S}: \R^{n} \rightrightarrows \R^{n}$ denotes the normal cone operator for the set $S \subseteq \R^{n}$, i.e., $\text{N}_{S}(x) = \varnothing$ if $x \notin S$, $\left\{v \in \mathbb{R}^{n} | \sup _{z \in S} v^{\top}(z-x) \leq 0\right\}$ otherwise.} $\operatorname{Id}$ is the identity operator; $I_n$ is the identity matrix of dimension $n$ and $ \bfs{0}_n$ is vector column of $n$ zeros; {their index is omitted where the dimensions can be deduced from context. The unit ball of appropriate dimensions depending on context is denoted with $\mathbb{B}$.} A continuous function $\gamma: \R_+ \rightarrow \R_+$ is of class $\mathcal{K}$ if it is zero at zero and strictly increasing. A continuous function $\alpha: \R_+ \rightarrow \R_+$ is of class $\mathcal{L}$ if is non-increasing and converges to zero as its arguments grows unbounded. A continuous function $\beta: \R_+ \times \R_+ \rightarrow \R_+$ is of class $\mathcal{KL}$ if it is of class $\mathcal{K}$ in the first argument and of class $\mathcal{L}$ in the second argument. UGES and SGPAS refer to uniform global exponential stability and semi-global practical asymptotic stability, respectively, as defined in \cite[Def. 2.2, 2.3]{poveda2021robust}.

\section{Discrete-time averaging}
We consider the following discrete-time system written in hybrid system notation \cite[Eq. 1.1, 1.2]{goebel2012hybrid}
\begin{align}
&\left\{\begin{array}{ll}
 {u}^+ &= u + \varepsilon G(u, \mu) \\
 {\mu}^+ &= M(u, \mu)
\end{array} \right., \begin{array}{rl}
 (u, \mu) &\in \mathcal{U} \times \Omega.
\end{array} \label{eq: non-averaged system}
\end{align} 

where $u$ and $\mu$ are the state variables, $\mathcal{U} \subset \R^{m}$, $\Omega \subset \R^l$, $G: \mathcal{U} \times \Omega \rightarrow \R^m$ and $M: \mathcal{U} \times \Omega \rightarrow \Omega$ are the state jump functions for the states $u$ and $\mu$ respectively, and $\varepsilon > 0$ is a small parameter. Furthermore, the mapping $G$ is parametrized by a small parameter $\gamma > 0$, i.e. $G = G_\gamma$, but for notational convenience, this dependence is not written explicitly in the equations. \\ \\

 Next we consider now the trajectories of the state $\mu$ oscillate indefinitely and in turn create oscillations in the trajectory of the state $u$. An equivalent system that produces trajectories $\bar{u}$ without oscillations should be easier to analyze. We refer to such systems as averaged systems \cite[Eq. 2.2.12]{bai1988averaging}, \cite[Eq. 8.33]{mareels1996adaptive}, and we focus on those of the following form:
\begin{align}
   \left\{ \begin{array}{ll}
         \tilde{u}^+ = \tilde{u} + \varepsilon G_{\textup{avg}}(\tilde{u}, \tilde{\mu}) \\
         {\tilde{\mu}}^+ = M(\tilde{u}, \tilde{\mu})
    \end{array} \right., \begin{array}{rl}
         (\tilde{u}, \tilde{\mu}) &\in \mathcal{U} \times \Omega,
    \end{array}    
    \label{eq: average system}
\end{align}
where $G_{\textup{avg}} : \mathcal{U} \times \Omega \rightarrow \R^m$ and is also parametrized by $\gamma > 0$. Unlike \cite[Thm. 2.2.4]{bai1988averaging}, \cite[Thm. 8.2.28]{mareels1996adaptive}, we take into consideration the case where the function $G_{\textup{avg}}$ depends on the fast state $\tilde{\mu}$, not only on $\tilde{u}$. \\ \\

To postulate the required relation between the function $G$ and the mapping $G_{\textup{avg}}$, we should introduce an auxiliary system that describes the behaviour of system \nref{eq: non-averaged system} when the state $u$ is kept constant, i.e. $\varepsilon = 0$, the so-called boundary layer system \cite[Eq. 6]{wang2012analysis}:

\begin{align}
&\left\{\begin{array}{l}
         {{u}_{\textup{bl}}^+} = u_{\textup{bl}} \\
         {{\mu}_{\textup{bl}}^+} = M(u_{\textup{bl}}, \mu_{\textup{bl}}) 
    \end{array} \right., (u_{\textup{bl}}, \mu_{\textup{bl}}) \in \mathcal{U} \times \Omega.
\end{align} \label{eq: boundary layer system}

Thus, a function $G_{\textup{avg}}$ is called an average of the mapping $G$ with the boundary layer dynamics in \nref{eq: boundary layer system} if the following condition holds true:
\begin{assum} \label{assum: boundness of the average aprox}
For any compact set $K \subset U$ and any solution $(u_{\textup{bl}}, \mu_{\textup{bl}})$ of \nref{eq: boundary layer system} where $u_{\textup{bl}}$ is contained in the compact set $K$, it holds that:
\begin{align}
    \n{\frac{1}{N}\sum_{i = 0}^{N - 1}\zs{G\z{u_{\textup{bl}}(i),\, \mu_{\textup{bl}}(i)} - G_{\textup{avg}}\z{u_{\textup{bl}}(i), \mu_{\textup{bl}}(i)}}} \red 
    \leq \sigma\z{N}, \label{eq: average bound}
\end{align} 
    for some function $\sigma: \R_+ \rightarrow \R_+$ of class $\mathcal{L}$. \kraj
\end{assum}
In plain words, in Assumption \ref{assum: boundness of the average aprox} we postulate that by using more samples over time, our approximation of the mapping $G$ becomes better. Furthermore, let us assume local Lipschitz continuity of the mappings as in \cite[Assum. 8.2.13]{mareels1996adaptive}, \cite[Eq. 1, Def. 1]{wang2009input} and compactness :
\begin{assum} \label{assum: continuity of jump and averag functions}
The functions $G$, $M$ and $G_{\textup{avg}}$  in \nref{eq: non-averaged system}, \nref{eq: average system} are continuous in their arguments and locally bounded; the mapping $G_{\textup{avg}}$ is locally Lipschitz continuous in its first argument. The set $\Omega$ is compact.\kraj
\end{assum}

The averaging method can be used in unison with other algorithms via time-scale separation. In such cases, often the averaged system does not exponentially or asymptotically stabilize the equilibrium as in \cite[Thm. 2.22]{bai1988averaging}, \cite[Thm. 8.2.28]{mareels1996adaptive},  due to the introduction of perturbations from other dynamics. Here, we assume the weaker property of semi-global practical stability of the set $\mathcal{A} \times \Omega$ under $v$-perturbed dynamics of the averaged system, where the perturbations are given by the dynamical system
\begin{align}
    v^+ = U(\tilde{u}, \tilde{\mu}, v) \label{eq: averaging additional dynamics}
\end{align}
with $v \in \R^m$, $U: \mathcal{U}\times\Omega\times\R^m \rightarrow \R^m$ being a function parametrized by a some $\varepsilon > 0$. 

\begin{assum}\label{assum: boundedness of hidden dynamics}
Consider the system in \nref{eq: average system} and perturbation dynamics in \nref{eq: averaging additional dynamics}, respectively. For any set $K \subset \mathcal{U}$, and all trajectories $(\tilde{u}, \tilde{\mu})$ contained in $K \times \Omega$, there exists a function of class $\mathcal{K}$, $\oln{v}$, such that $\max_{k \in \dom(v)}\n{v(k)} \leq \oln{v}(\varepsilon)$.  \kraj
\end{assum}

\begin{assum}\label{assum: stability of the average system}
The set $\mathcal{A} \times \Omega$ is SGPAS as $\gamma \rightarrow 0$ for the dynamics in \nref{eq: average system}, perturbations in \nref{eq: averaging additional dynamics}, and the corresponding Lyapunov function $V_{\textup{a}}$ satisfies:

\begin{subequations}
\begin{align}
    &\underline{\alpha}_{\textup{a}}\z{\n{z}_{\mathcal{A}}} \leq V_{\textup{a}}(z, \mu) \leq  \oln{\alpha}_{\textup{a}}\z{\n{z}_{ \mathcal{A}}}  \label{eq: assum lyapunov bound} \\
    &V_{\textup{a}}(z^+, \mu^+) - V_{\textup{a}}(z, \mu) \leq - \tilde{\alpha}_{\varepsilon}\z{{\varepsilon}}\alpha_{\textup{a}}\z{\n{z}_{\mathcal{A}}} \red
    &\text{ for } \n{z}_\mathcal{A} \geq \alpha_{\gamma}(\gamma),
\end{align}
\end{subequations} 
where $z = \tilde{u} + v$, $\underline{\alpha}_{\textup{a}}, \overline{\alpha}_{\textup{a}}, \tilde{\alpha}_{\varepsilon}, \alpha_{\textup{a}}, \alpha_{\gamma}$ are functions of class $\mathcal{K}$, and the function $\frac{\varepsilon}{\tilde{\alpha}_{\varepsilon}\z{\varepsilon}}$ is bounded for $\varepsilon \in (0, \oln{\varepsilon})$.\kraj
\end{assum}                          

Under these assumptions, we claim that the original system is semi-global practically asymptotically stable, as formalized next:
\begin{theorem}\label{thm: averaging theorem}
Let Assumptions \ref{assum: boundness of the average aprox}, \ref{assum: continuity of jump and averag functions}, \ref{assum: boundedness of hidden dynamics} and \ref{assum: stability of the average system} hold. The set $\mathcal{A}\times\Omega$ is SGPAS as $(\varepsilon, \gamma)\rightarrow 0$ for the discrete dynamics in \nref{eq: non-averaged system} with perturbations in \nref{eq: averaging additional dynamics}. The corresponding Lyapunov function $V_{\textup{a}}$ satisfies:
\begin{align*}
    &\underline{\alpha}_{\textup{a}}\z{\n{\xi}_{\mathcal{A}}} \leq V_{\textup{a}}(\xi, \mu) \leq  \oln{\alpha}_{\textup{a}}\z{\n{\xi}_{\mathcal{A}}}  \\
    &V_{\textup{a}}(\xi^+, \mu^+) - V_{\textup{a}}(\xi, \mu) \leq - \hat{\alpha}_{\varepsilon}\z{{\varepsilon}}\alpha_{\textup{a}}\z{\n{\xi}_{\mathcal{A}}} \\
    &\text{ for } \n{\xi}_{\mathcal{A}} \geq \max\{\alpha_{\gamma}(\gamma), {\alpha}_{\varepsilon}(\varepsilon)\},
\end{align*}
where $\xi \coloneqq u + v + \eta$, $\eta$ is the perturbation state with dynamics
\begin{align}
&\eta^+ = (1 - \varepsilon)\eta + \varepsilon [G_{\textup{avg}}(u, \mu) - G(u, \mu)], \red
&\max_{k \in \N}\n{\eta(k)} \leq \oln{\eta}(\varepsilon),
\end{align}
the $v$ dynamics are given by \nref{eq: averaging additional dynamics}, and $\hat{\alpha}_{\varepsilon}$, ${\alpha}_{\varepsilon}$, $\oln{\eta}$ are functions of class $\mathcal{K}$.
\kraj
\end{theorem}

\begin{pf} 
See Appendix \ref{proof: averaging theorem}. \krajdokaz
\end{pf}

\section{Applications of the averaging theorem}
In this section, we apply our averaging theorem, Theorem \ref{thm: averaging theorem}, to derive two novel convergence results for NEPs. First, we propose a zeroth-order algorithm for solving strongly monotone NEPs with local constraints \emph{in discrete time}. Secondly, we propose an algorithm for solving strongly monotone unconstrained NEPs where the agents sample their states \emph{asynchronously}.

\subsection{Zeroth-order discrete time forward-backward algorithm}\label{average example 1}
Let us consider a multi-agent system with $M$ agents indexed by $i \in \mathcal{I} \coloneqq \{1, 2, \dots M\}$, each with cost function 
\begin{align}
    J_i(x_i, \boldsymbol{x}_{-i}), \label{eq: problem1 def}
\end{align}

where $x_i \in \Omega_i \subset \R^{m_i}$ is the decision variable,  $J_i: \R^{m_i} \times \R^{m_{-i}} \rightarrow \R$, $m \coloneqq \sum_{j \in \mathcal{I}} m_j$, $m_{-i} \coloneqq \sum_{j \neq i} m_j$, $\Omega \coloneqq \Omega_i \times \dots \times \Omega_N$. Formally, let the goal of each agent be to reach a steady state that minimizes their own cost function, i.e.,
\begin{align}
\forall i \in \mathcal{I}:\ &\min_{x_i \in \Omega_i}  J_i(x_i, \boldsymbol{x}_{-i}). \label{def: dyn_game} \end{align}
A popular solution to this problem is the so-called Nash equilibrium:
\begin{defn}[Nash equilibrium]\hfill\quad
A set of decision variables $\bfs{x}^*\coloneqq\col{x_i^*}_{i \in \mathcal{I}}$ is a Nash equilibrium if, for all $i \in \mathcal{I}$,
\begin{align}
    x_{i}^{*} \in \underset{v_{i} \in \Omega_i}{\operatorname{argmin}}\ J_{i}\left(v_{i}, \bfs{x}_{-i}^{*}\right). \kraj\label{def: ne}
\end{align} \end{defn}
A fundamental mapping in NEPs is the pseudogradient mapping $F: \R^m \rightarrow \R^m$, which is defined as:
\begin{align}
    F(\boldsymbol{x}):=\operatorname{col}\left(\left(\nabla_{x_{i}} J_{i}\left(x_{i}, \bfs{x}_{-i}\right)\right)_{i \in \mathcal{I}}\right). \label{eq: pseudogradient}
\end{align}
Let us also define $C_{\textup{F}} \coloneqq \oln{\operatorname{co}} \{ F(\Omega) \}$, the convex hull of the image of the pseudogradient. To ensure existence and uniqueness of the Nash equilibrium, we  assume certain regularity properties \cite[Thm. 4.3]{bacsar1998dynamic}:
\begin{assum}\label{assum: regularity of the problem}
For each $i \in \mathcal{I}$, the function $J_i$ in \nref{eq: problem1 def} is continuously differentiable in $x_i$ and continuous in $\bfs{x}_{-i}$; the function $J_{i}\left(\cdot, \bfs{x}_{-i}\right)$ is strictly convex for every fixed $\bfs{x}_{-i}$.\kraj
\end{assum}
Furthermore, let us assume that no agent can compute their part of the the pseudogradient $F$ directly, but they can only measure their instantaneous cost  $h_i = J_i(x_i, \bfs{x}_{-i})$, a common assumption in extremum-seeking problems \cite{krilavsevic2021learning}, \cite{poveda2020fixed}, \cite{poveda2017framework}, \cite{stankovic2011distributed}. The full-information problem where $F$ is known can be solved in many ways, depending on the technical assumptions on the problem data. Here we choose to study a simple forward-backward algorithm \cite[Equ. 26.14]{bauschke2011convex}:
\begin{align}
    \bfs{x}^+ = (1 - \lambda)\bfs{x} + \lambda \proj_{C}\z{\bfs{x} - \gamma F(\bfs{x})}, \label{eq: example1 fb algorithm}
\end{align}
for which the Lyapunov function $V(\bfs{x}) = \n{\bfs{x} - \bfs{x}^*}^2$ satisfies the inequality
\begin{align}
    &V(\bfs{x}^+) - V(\bfs{x}) \leq - \lambda(1 - c)(2 - \lambda c)V(\bfs{x}), \label{eq: example1 algorithm lyapunov}
\end{align}
where $c \coloneqq \frac{\sqrt{1 + \gamma^2 L^2}}{1 + \gamma \mu_F}$ and $\bfs{x}^*$ is the Nash equilibrium of the game in \nref{eq: problem1 def}. We note that this Lyapunov function satisfies Assumption \ref{assum: stability of the average system}.\\
A naïve approach to adapting the algorithm in \nref{eq: example1 fb algorithm} for zeroth-order implementation would be to use a gradient estimation scheme as in \cite{choi2002extremum}, \cite{stankovic2011distributed} and plug in the estimate directly into \nref{eq: example1 fb algorithm}. However, because of the projection, Assumption \ref{assum: boundness of the average aprox} would not be satisfied. Thus, an additional time-scale separation is hereby proposed:

\begin{align}
&\left\{ \begin{array}{ll}
              \bfs{x}^+ &= (1 - \alpha\beta)\bfs{x} + \alpha\beta \proj_{C}\z{\bfs{x} - \gamma \bfs{\xi}} \\
            \bfs{\xi}^+ &= (1 - \alpha)\bfs{\xi} + \alpha 2A^{-1} J(\bfs{x} + A\mathbb{D}\bfs{\mu})\mathbb{D}\bfs{\mu} \\
            \bfs{\mu}^+ &= \mathcal{R}\bfs{\mu}
        \end{array}\right.,\label{eq: example1 algorithm}
\end{align}

where $\bfs{\xi} \in \R^m$ are filter states, $\bfs{\mu} \in \mathbb{S}^{m}$  are the oscillator states, $\alpha, \beta > 0$ are small time-scale separation parameters, $\mathcal{R} \coloneqq \Diag{(\mathcal{R}_i)_{i \in \mathcal{I}}}$, $\mathcal{R}_i \coloneqq \Diag{\sm{\cos(\omega_i^j) & -\sin(\omega_i^j) \\ \sin(\omega_i^j) & \cos(\omega_i^j)}_{j \leq m_i}}$, $\omega_i^j > 0$ for all $i$ and $j$, $\mathbb{D} \in \R^{m \times 2m}$ is a matrix that selects every odd row from the vector of size $2m$, $a_i > 0$ are small perturbation amplitude parameters, $A \coloneqq \diag{(a_i)_{i \leq m}}$ and $J(\bfs{x}) = \Diag{(J_i(x_i, \bfs{x}_{-i})I_{m_i})_{i \in \mathcal{I}}}$. We claim that the dynamics in \nref{eq: example1 algorithm} render the set $\{\bfs{x}^*\} \times C_{\textup{F}} \times \mathbb{S}^m$ practically stable. To the best of our knowledge, it is not possible to prove convergence of the algorithm in \nref{eq: example1 algorithm}, using the current averaging theory for discrete-time systems, since \cite[Thm. 8.2.28]{mareels1996adaptive}, \cite[Thm. 2.22]{bai1988averaging} require exponential stability of the origin via the averaged system, and \cite[Thm. 2]{wang2009input} does not incorporate boundary-layer dynamics. We claim that under the strong monotonicity assumption of the pseudogradient, and a proper choice of the perturbation frequencies, the algorithm in \nref{eq: example1 algorithm} converges to a Nash equilibrium.

\begin{assum}\label{assum: exmaple2 strong monotinicty}
    The pseudogradient mapping $F$ is $\mu_\text{f}$-strongly monotone and $L$-Lipschitz continuous, i.e. $\vprod{\bfs{x} - \bfs{y}}{F(\bfs{x}) - F(\bfs{y})} \geq \mu_\text{f}\n{\bfs{x} -\bfs{y}}$, $\n{F(\bfs{x
}) - F(\bfs{y})} \leq L \n{\bfs{x} - \bfs{y}}$, for all $(\bfs{x}, \bfs{y}) \in \R^{2m}$.\kraj
\end{assum}
\begin{assum}\label{assum: exmaple1 set convexity}
    The sets $(\Omega_i)_{i \in \mathcal{I}}$ are convex, closed and bounded.\kraj
\end{assum}

\begin{assum}\label{assum: example 1 different freq}
    The rotational frequencies of each agent $i$, $\bfs{\omega}_i = \col{(\omega_i^j)_{j \leq m_i}}$, are chosen so that $\omega_i^j \pm \omega_u^v \neq 2\pi z', z' \in \Z$, for every $u \in \mathcal{I}$, for every $j \in \{1, \dots, m_i\}$, for every $v \in \{1, \dots, m_u\}$, apart for the case when $i = u$ and $j = v$. \kraj
\end{assum}

\begin{theorem}\label{theorem: example1 convergence}
Let Assumptions \ref{assum: regularity of the problem}, \ref{assum: exmaple2 strong monotinicty}, \ref{assum: exmaple1 set convexity} and \ref{assum: example 1 different freq} hold. The set $\{\bfs{x}^*\} \times C_{\textup{F}} \times \mathbb{S}^m$ is SGPAS as $(\alpha, \oln{a}, \beta) \rightarrow 0$ for the dynamics in \nref{eq: example1 algorithm}.\kraj
\end{theorem}

\begin{pf}
See Appendix \ref{proof: thm convergence full info}. \krajdokaz
\end{pf}

\subsection{Asynchronous zeroth-order discrete time forward algorithm}
\subsubsection{First-order information feedback}
We now consider the same NEP as in Section \ref{average example 1} with $\Omega_i \coloneqq \R^{m_i}$, but where the agent are \emph{asynchronous}, i.e. each agent samples their states independently of others without a global clock for synchronization. For ease of exposition, we assume that the initial conditions are chosen so that simultaneous sampling never occurs. In the full-information case, such algorithm can be represented in the following form:
\begin{subequations}\label{eq: example 2 hybrid rep single agent full info}
\begin{align}
    &\col{\dot{x}_i, \dot{\tau}_i, \dot{\kappa_i}, \dot{t}} = \col{\bfs{0
    }, \tfrac{1}{T_i}, 0, 1} \red 
    &\text{ if } \col{x_i, \tau_i, \kappa_i, t}\in \R^{m_i} \times [0, 1] \times \N \times \R \\
    &\left\{ \begin{array}{ll}
    x_i^{+} &= x_i - \, \alpha \nabla_{x_i}J_i(x_i, \bfs{x}_{-i}) \\
    \tau_i^+ &= 0 \\ 
    \kappa_i^+ &= \kappa_i + 1 \\
    t^+ &= t
         \end{array}\right.
 \red 
 &\text{ if } \col{x_i, \tau_i, \kappa_i, t}\in \R^{m_i} \times \{1\}\times \N \times \R ,
\end{align}
\end{subequations}

which in collective form reads as 
\begin{subequations}\label{eq: example2 hybrid rep full info}
\begin{align} 
    &\col{\dot{\bfs{x}}, \dot{\bfs{\tau}}, \dot{\bfs{\kappa}}, \dot{k}, \dot{t}} = \col{\bfs{0}, \bfs{T}^{-1}, \bfs{0}, 0, 1} \red 
    &\text{ if } \col{\bfs{x}, \bfs{\tau}, \bfs{\kappa}, k, t}\in \R^{m} \times \mathcal{T} \times \N^{N + 1} \times \R \\
    &\left\{ \begin{array}{ll}
        \bfs{x}^{+} &= \bfs{x} - \alpha S_{x}(\bfs{\tau})  F(\bfs{x})\\
        \bfs{\tau}^+ &= (I - S_{\tau}(\bfs{\tau}))\bfs{\tau} \\
        \bfs{\kappa}^+ &= \bfs{\kappa} + S_\tau(\bfs{\tau})\\
         k^+ &= k + 1 \\ 
         t^+ &= t\end{array}\right.\red 
    &\text{ if } \col{\bfs{x}, \bfs{\tau}, \bfs{\kappa}, k, t}\in \R^{m} \times \mathcal{T}_\textup{R} \times \N^{N + 1} \times \R,
\end{align}
\end{subequations}

where $\tau_i$ are timer states, $t$ is the ``experienced" global time, $\bfs{T}^{-1} \coloneqq \col{(\tfrac{1}{T_i})_{i \in \mathcal{I}}}$ is the vector of inverse sampling times, $\mathcal{T} \subset [0, 1]^N$ is a closed invariant set in which all of the timers evolve and it excludes the initial conditions and their neighborhood for which we have concurrent sampling, $\mathcal{T}_\textup{R} \coloneqq \z{\cup_{i\in\mathcal{I}} [0, 1]^{i - 1} \times \{1\} \times [0, 1]^{N - i}}\cap\mathcal{T}$ is the set of timer intervals where at least one agent has triggered its sampling, $\kappa_i$ are private event counters, $k$ is the global event counter, $S_{x}: \mathcal{T} \rightarrow \R^{m \times m}$ and $S_{\tau}: \mathcal{T} \rightarrow \R^{N \times N}$ are continuous functions that output diagonal matrices with ones on the positions that correspond to states and timers of agents with $\tau_i = 1$, respectively, while other elements are equal to zero, when evaluating at $\bfs{\tau} \in \mathcal{T}_\textup{R}$. We note that the functions $S_{x}, S_{\tau}$ are introduced only to write down the algorithm in the collective form, while the agents themselves do not require them for their dynamics and in fact just follow  \nref{eq: example 2 hybrid rep single agent full info}. Furthermore, the counter states $\kappa_i, k$ and global time $t$ are not necessary for the algorithm convergence, yet they help with understanding the setup of the algorithm. We choose to represent the algorithm in the hybrid dynamical system framework to replicate the behaviour of sampled systems with different sampling periods, and to see its effects on the functions $S_x, S_\tau$. Later, we represent and model the system as a fully discrete-time system in order to study convergence.\\ First, we show that the solution $\bfs{\tau}(t, j)$ is periodic and that implies that $S_{x}(\bfs{\tau}(t, j))$ and $S_{\tau}(\bfs{\tau}(t, j))$ are also periodic. We make the following assumption:
\begin{assum}\label{assumption example2 sampling ratio}
    There exist natural numbers $\z{p_i}_{i \in \mathcal{I}}$, such that the proportion $T_1 : T_2 : \dots : T_N = p_1 : p_2 \dots : p_N$ holds, where $\z{T_i}_{i \in \mathcal{I}}$ are the sampling times. \kraj
\end{assum}

\begin{lemma}\label{lemma: example 2 periodicity}
Let Assumption \ref{assumption example2 sampling ratio} hold. Denote $r_i = \tfrac{p}{p_i}$ and $r = \sum_{i \in \mathcal{I}} r_i$, where $p$ is the least common multiple of $(p_i)_{i \in \mathcal{I}}$. For any trajectory  $S_{x}(\bfs{\tau}(t, j))$ and $S_{\tau}(\bfs{\tau}(t, j))$, where $\bfs{\tau}(t, j)$ is a solution of the system in \nref{eq: example2 hybrid rep full info}, there exists $T > 0$ such that $S_{x}(\bfs{\tau}(t, j)) = S_{x}(\bfs{\tau}(t + T, j + r))$ and $S_{\tau}(\bfs{\tau}(t, j)) = S_{\tau}(\bfs{\tau}(t + T, j + r))$ for all $(t, j) \in \dom(\bfs{\tau})$ such that a jump occurred at time $t$.\kraj
\end{lemma}
\begin{pf}
    See Appendix \ref{appendix: proof of (T,q) period}. \krajdokaz
\end{pf}
Because the values of $S_x$ and $S_\tau$ are used only during jumps, we define 
\begin{align}
    \hat{S}_x(k; \bfs{\tau}(0, 0)) = S_x(\max_{(t\in \dom(\bfs{\tau}(\cdot, k))} t, k)\\
    \hat{S}_{\tau}(k; \bfs{\tau}(0, 0)) = S_{\tau}(\max_{(t\in \dom(\bfs{\tau}(\cdot, k))} t, k),
\end{align}
where functions $\hat{S}_x: \N \rightarrow \R^{m \times m}$ and $\hat{S}_\tau: \N \rightarrow \R^{N \times N}$ are parametrized by the vector of initial conditions of the timers, since different initial conditions can change the order in which the agents are sampling their actions. Due to Lemma, \ref{lemma: example 2 periodicity}, for every initial condition $\bfs{\tau}(0, 0) = \bfs{\tau}^0$, it follows that $\hat{S}_x(k, \bfs{\tau}^0) = \hat{S}_x(k + r, \bfs{\tau}^0)$ for all $k \in \N$.\\

Now we consider the following discrete time systems

\begin{align}
    x_i(k + 1) &= x_i(k) - \, \alpha \hat{S}^i_x(k; \bfs{\tau}_0)\nabla_{x_i}J_i(x_i(k), \bfs{x}_{-i}(k))
\end{align}

which in collective form read as 
\begin{align} 
    \bfs{x}(k + 1)&= \bfs{x}(k) - \alpha \hat{S}_x\z{k; \bfs{\tau}_0}  F(\bfs{x}(k)), \label{eq: example2 discrete-time collective full info}
\end{align}

where the function $\hat{S}^i : \N \rightarrow \R^{m_i \times m_i}$ returns the rows of $\hat{S}_{x}^i(k; \bfs{\tau}_0)$ corresponding to agent $i$. We can show that for every solution of \nref{eq: example2 discrete-time collective full info} there exists a corresponding solution of \nref{eq: example2 hybrid rep full info} and vice versa. We claim that under the strong monotonicity assumption, an additional regularity assumption due to the unboundedness of the decision set, and proper choice of the parameter $\alpha$, the dynamics in \nref{eq: example2 discrete-time collective full info} converge to the solution of the game, with the same minimal convergence rate, regardless of the initial conditions of the timers.

\begin{assum}\label{assum: regularity of the problem, unbounded}
For each $i \in \mathcal{I}$, the function $J_i(\cdot, \bfs{x}_{-i})$ in \nref{eq: problem1 def} is radially unbounded for every fixed $\bfs{x}_{-i}$.\kraj
\end{assum}

\begin{theorem}\label{thm: essentisally cyclic forward algorithm}
Let Assumptions \ref{assum: regularity of the problem}, \ref{assum: exmaple2 strong monotinicty}, \ref{assumption example2 sampling ratio} and \ref{assum: regularity of the problem, unbounded} hold. Then, for all vectors of initial conditions $\bfs{\tau}_0$, there exists $\alpha^*$, such that for any $\alpha \in (0, \alpha^*)$, the NE solution $\bfs{x}^*$ is UGES for the dynamics in \nref{eq: example2 discrete-time collective full info}. Furthermore, the corresponding Lyapunov function satisfies Assumption \ref{assum: stability of the average system}. \kraj
\end{theorem}
 
\begin{pf}
See Appendix \ref{appendix: proof essentisally cyclic forward algorithm}.\krajdokaz
\end{pf}

Moreover, for the hybrid system representation in \nref{eq: example2 hybrid rep full info}, since the trajectories of $(\bfs{\tau}, \bfs{\kappa}, k, t)$ are invariant to the set $\mathcal{T} \times \N^{N + 1} \times \R$, and by the structure of the flow and jump sets in \nref{eq: example2 hybrid rep full info} that assures complete solutions with unbounded time and jump domains, it follows that the dynamics in \nref{eq: example2 hybrid rep full info} render the set $\{\bfs{x}^*\} \times \mathcal{T} \times \N^{N + 1} \times \R$ UGES, as formalized next.

\begin{corollary}\label{corollary: stability of full info hybrid asynchronus}
Let the Assumptions \ref{assum: regularity of the problem}, \ref{assum: exmaple2 strong monotinicty}, \ref{assumption example2 sampling ratio} and \ref{assum: regularity of the problem, unbounded} hold. Then, the set $\{\bfs{x}^*\} \times \mathcal{T} \times \N^{N + 1} \times \R$ is UGES for the dynamics in \nref{eq: example2 hybrid rep full info}. Furthermore, the corresponding Lyapunov function satisfies Assumption \ref{assum: stability of the average system}. \kraj
\end{corollary}

\subsubsection{Zeroth-order information feedback}\label{sec: zeroth-order asynch pseudogradient descent}
Now, consider that each agent only has access to the measurements of the cost function. They can modify the algorithm in \nref{eq: example 2 hybrid rep single agent full info} by implementing a pseudogradient estimation scheme, similar to the one in Equation \nref{eq: example1 algorithm}:

\begin{subequations} 
\begin{align}
    &\col{\dot{x}_i, \dot{\xi}_i , \dot{\mu}_i, \dot{\tau}_i, \dot{\kappa_i}, \dot{t}} = \col{\bfs{0}, \bfs{0}, \bfs{0}, \tfrac{1}{T_i}, 0, 1} \\
    &\text{ if } \col{x_i, \xi_i, \mu_i, \tau_i, \kappa_i, t}\in \R^{2m_i}  \times \mathbb{S}^m \times [0, 1] \times \N \times \R, \red
    &\left\{ \begin{array}{ll}
        x_i^{+} = x_i - \, \alpha\beta \xi_i \\
        \xi_i^+ = (1 - \alpha)\xi + \alpha \tfrac{2}{a_i} J_i(x + A\mathbb{D}\mu)\mathbb{D}_i\mu_i \\ 
        \mu_i^+ = \mathcal{R}_i\mu_i \\
        \tau_i^+ =  0 \\ 
        \kappa_i^+ = \kappa_i + 1 \\
        t^+ = t \end{array}\right.\\
    & \text{ if } \col{x_i, \xi_i, \mu_i, \tau_i, \kappa_i, t}\in \R^{2m_i} \times \mathbb{S}^m \times \{1\}\times \N \times \R,\nonumber
\end{align}
\end{subequations}

which in the collective form reads as:

\begin{subequations} \label{eq: example2 algorithm}
\begin{align} 
    &\operatorname{col}(\dot{\bfs{x}}, \dot{\bfs{\xi}}, \dot{\bfs{\mu}}, \dot{\bfs{\tau}}, \dot{\bfs{\kappa}}, \dot{k}, \dot{t}) = \col{\bfs{0}, \bfs{0}, \bfs{0}, \bfs{T}^{-1}, \bfs{0}, 0, 1} \\
    &\text{ if } \col{\bfs{x}, \bfs{\xi}, \bfs{\mu}, \bfs{\tau}, \bfs{\kappa}, k , t} \in \R^{2m} \times \mathbb{S}^{m} \times \mathcal{T} \times \N^{N + 1}\times \R, \red
    &\left\{ \begin{array}{ll}
        \bfs{x}^{+} = \bfs{x} - \alpha \beta S_x(\bfs{\tau})\bfs{\xi}\\ 
        \bfs{\xi}^{+} = \bfs{\xi} + \alpha S_x(\bfs{\tau})\z{2A^{-1} J(\bfs{x} + A\mathbb{D}\bfs{\mu})\mathbb{D}\bfs{\mu} - \bfs{\xi}} \\
        \bfs{\mu}^{+} = (I - S_{\mu}(\bfs{\tau}))\bfs{\mu} + S_{\mu}(\bfs{\tau}))\mathcal{R}\bfs{\mu} \\
        \bfs{\tau}^+ = (I - S_{\tau}(\bfs{\tau}))\bfs{\tau} \\ 
        \bfs{\kappa}^+ = \bfs{\kappa} + S_\tau(\bfs{\tau})\\
        k^+ = k + 1 \\
        t^+ = t
             \end{array}\right. \\
    &\text{ if } \col{\bfs{x}, \bfs{\xi}, \bfs{\mu}, \bfs{\tau}, \bfs{\kappa}, k, t} \in \R^{2m} \times \mathbb{S}^{m}\times \mathcal{T}_\textup{R} \times \N^{N + 1} \times \R,\nonumber
\end{align} 
\end{subequations}

where $S_{\mu}: \mathcal{T}\rightarrow \R^{2m \times 2m}$ is a continuous functions that outputs a diagonal matrix with ones on the positions that correspond to oscillator states of agents with $\tau_i = 1$, while other elements are equal to zero, when evaluating at $\bfs{\tau} \in \mathcal{T}_\textup{R}$, and other notation is defined as in \nref{eq: example1 algorithm} and \nref{eq: example2 hybrid rep full info}.\\

Under the assumption of properly chosen perturbation frequencies, we claim semi-global practical stability of the set of solutions.
\begin{assum}\label{assum: example 2 different freq}
    The rotational frequencies of every agent $i$, $\bfs{\omega}_i = \col{(\omega_i^j)_{j \leq m_i}}$, are chosen so that $\omega_i^j r_i \pm \omega_u^v r_j \neq 2\pi z', z' \in \Z$, $r_i = \tfrac{p}{p_i}, r_j = \tfrac{p}{p_j}$, for every $u \in \mathcal{I}$, for every $j \in \{1, \dots, m_i\}$, for every $v \in \{1, \dots, m_u\}$, apart for the case when $i = u$ and $j = v$. \kraj
\end{assum}
\begin{theorem}\label{theorem: example2 convergence}
Let the Assumptions \ref{assum: regularity of the problem}, \ref{assum: exmaple2 strong monotinicty}, \ref{assumption example2 sampling ratio}, \ref{assum: example 2 different freq} hold. The set $\{x^*\} \times \R^m \times \mathbb{S}^m \times \mathcal{T} \times \N^{N + 1} \times \R$ is SGPAS as $(\alpha, \oln{a}, \beta) \rightarrow 0$ for the dynamics in \nref{eq: example2 algorithm}.\kraj
\end{theorem}

\begin{pf}
The result is proven by following the same steps as the proof of Theorem \ref{theorem: example1 convergence} and by using system in \nref{eq: example2 hybrid rep full info} with additional filtering state $\bfs{\xi}$ like in \nref{eq: example1 average system2} as the second averaged system. \krajdokaz
\end{pf}

\section{Illustrative example}

The connectivity control problem has been considered in \cite{stankovic2011distributed} as a Nash equilibrium problem. In many practical scenarios, multi-agent systems, besides their primary objective, are designed to uphold certain connectivity as their secondary objective.  In what follows, we consider a similar problem in which each agent is tasked with finding a source of an unknown signal while maintaining certain connectivity. Unlike \cite{stankovic2011distributed}, we only consider the case without vehicle dynamics. \\
Consider a system consisting of multiple agents indexed by $i \in \mathcal{I} \coloneqq\{1, \dots N\}$. Each agent is tasked with locating a source of a unique unknown signal. The strength of all signals abides by the inverse-square law, i.e. proportional to $1/r^2$. Therefore, the inverse of the signal strength can be used as a cost function. Additionally, the agents must not drift apart from each other too much, as they should provide quick assistance to each other in case of critical failure. This is enforced by incorporating the signal strength of the fellows agents in the cost functions. Thus, we design the cost 
\begin{align}
    \forall i \in \mathcal{I}: 
    J_i(\bfs{x}) = \|x_i - x_i^s \|^2 + c\sum_{j \in \mathcal{I}_{-i}}\| x_i - x_j\|^2,
\end{align}
where $\mathcal{I}_{-i} \coloneqq \mathcal{I}\setminus \{i\}$, $c, b > 0$ and $x_i^s$ represents the position of the source assigned to agent $i$. Goal of each agent is to minimize their cost function, and the solution to this problem is a Nash equilibrium. Furthermore, agents are mutually independent so their sampling time are not synchronised. To solve this problem, we use the asynchronous pseudogradient descent algorithm in \nref{eq: example2 algorithm}.\\

\begin{figure}
    \centering
    \includegraphics{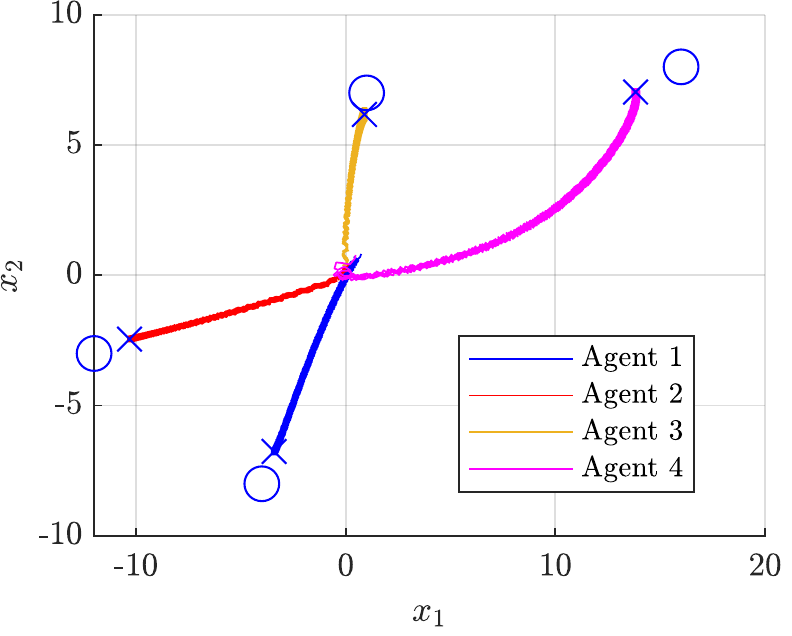}
    \caption{State trajectories in the $x_1-x_2$ plane. Circle symbols represent locations of the sources, while the $\times$ symbols represent locations of the NE. Perturbations signals are added to the states.}
    \label{fig:application x phase plane}
\end{figure}

\begin{figure}
    \centering
    \includegraphics{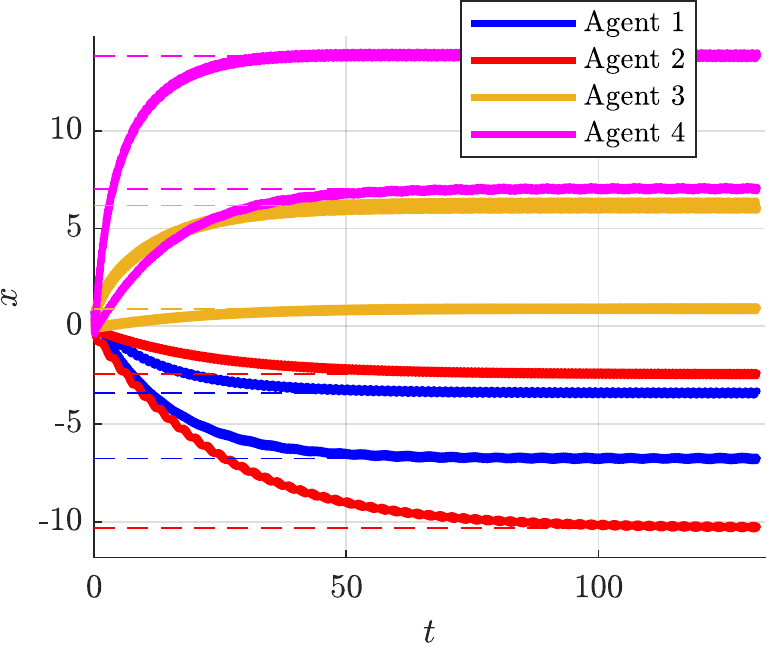}
    \caption{Time response of the states. The dashed lines correspond to the states of the Nash equilibrium.}
    \label{fig:application x time plot}
\end{figure}

For our numerical simulations, we choose the parameters: $x^s_1 = (-4, -8)$, $x^s_2 = (-12, -3)$, $x^s_3 = (1, 7)$, $x^s_4 = (16, 8)$, $c = 0.04$, $\gamma = 0.1$, $\alpha = 0.1$, $\beta = 0.003$, $a_i = 0$ for all $i$, $T = (0.01, 0.015, 0.02, 0.01)$, $\bfs{\tau}(0, 0) = (0, 0.002, 0.004, 0.006)$, the perturbation frequencies $\omega_i^j$ were chosen as different natural numbers with added random numbers of maximal amplitude of 0.5.\\
The numerical results are illustrated on Figures \ref{fig:application x phase plane} and \ref{fig:application x time plot}. We note that the trajectories converge to a small neighborhood of the Nash equilibrium. This can be partially attributed to the robustness properties of the pseudogradient descent with strongly monotone operators.


\section{Conclusion}
Averaging theory can be adapted for use in discrete systems with multiple timescales. Furthermore, strongly monotone Nash equilibrium problem with constrained action sets, or with asynchronous action sampling, can be solved via zeroth-order discrete-time algorithms that leverage novel averaging theory results.

\endlinechar=13
\bibliographystyle{plain}        
\bibliography{biblioteka}  

\endlinechar=-1
\appendix

\section{Proof of Theorem \ref{thm: averaging theorem}}\label{proof: averaging theorem}
\emph{Sketch of the proof:}  First, we show that under a change of coordinates, the system in \nref{eq: non-averaged system} can be represented as an inflated version of the averaged system in \nref{eq: average system}. Then we show that the inflation can be arbitrarily small for small enough $\varepsilon$. Finally, we use the stability properties of the averaged system and the bounded inflation property to prove SGPAS.\\ \\

By introducing an additional state, we construct the augmented system:

\begin{align}
&\left\{\begin{array}{l}
         {u}^+ = u + \varepsilon G(u, \mu) \\
         {\mu}^+ = M(u, \mu) \\
         \eta^+ = (1 - \varepsilon)\eta + \varepsilon\left[ G_{\textup{avg}}(u, \mu) - G(u, \mu) \right]\\
    \end{array}\right.\label{eq: augmented system}\\
&(u, \mu, \eta) \in \mathcal{U}  \times \Omega \times \R^m.\nonumber
\end{align} 
With a change of coordinates $u = \tilde{u} - \eta$, $\mu = \tilde{\mu}$ the system is transformed to:

\begin{align}
   &\left\{ \begin{array}{l}
         \tilde{u}^+ = \tilde{u} + \varepsilon G_{\textup{avg}}(\tilde{u} - \eta, \tilde{\mu}) - \varepsilon \eta\\
         \tilde{\mu}^+ = M(\tilde{u} - \eta, \tilde{\mu}) \\
         \eta^+ = (1 - \varepsilon)\eta + \varepsilon\left[G_{\textup{avg}}(\tilde{u} - \eta, \tilde{\mu}) - G(\tilde{u} - \eta, \tilde{\mu}) \right]
    \end{array}\right. \label{eq: transformed system} \\
    &(\tilde{u} - \eta, \tilde{\mu}, \eta) \in \mathcal{U} \times \Omega \times \R^m\nonumber
\end{align} 

We note that the $\tilde{u}$ dynamics in \nref{eq: transformed system} are perturbed dynamics of the averaged system in \nref{eq: average system} and $\n{\eta}$ is the upper bound on the perturbation amplitude. To prove our desired stability, we characterize the bound of this amplitude:

\begin{lemma}\label{lemma: boundedness of eta}
For every $a>0$ and compact set $K \in \mathcal{U}$, there exists $\varepsilon^*$ such that $\n{\eta(t, j)} < a$ holds for any  $\varepsilon \in (0, \varepsilon^*]$ and any trajectory of system in \nref{eq: transformed system} where $\tilde{u}$ is contained in the set $K$. \kraj
\end{lemma}

\begin{pf}
For the purposes of the proof, we construct the concatenated trajectory ($u^L,\, \mu^L$), which is created by taking solutions of length $L$ of the system in \nref{eq: boundary layer system} and concatenating them together. \\
We derive a similar bound to \nref{eq: average bound} for the concatenated trajectory type using Assumption \ref{assum: boundness of the average aprox}:
\begin{align}
    &\n{\sum_{i = 1}^{N}\zs{G\z{u^{L}_i, \mu_i^{L}} - G_{\textup{avg}}\z{u^{L}_i, \mu_i^{L}}}}  \red
    &\leq  \sum_{j = 1}^n \n{\sum_{i = 1 + L(j - 1)}^{1 + jL}\zs{G\z{u_i^{L}, \mu_i^{L}} - G_{\textup{avg}}\z{u^{L}_i, \mu_i^{L}}}} \red
    & + \n{\sum_{i = nL + 1}^N\zs{G\z{u_i^{L}, \mu_i^{L}} - G_{\textup{avg}}\z{u^{L}_i, \mu_i^{L}}}}\red
    & \leq nL \sigma(L) + (N - nL)\sigma(N - nL) \red
    & \leq N\sigma(L) + L\sigma(0) \label{eq: concat average bound}
\end{align}
Note that the bounds in \nref{eq: average bound} and \nref{eq: concat average bound} use (concatenated) boundary layer trajectories instead of the $u$ trajectory in \nref{eq: non-averaged system}. In order to use the bound in \nref{eq: concat average bound}, we rewrite the $\eta$ dynamics in \nref{eq: augmented system} as
\begin{align}
    \eta^+ = (1 - \varepsilon)\eta + v_1 + v_2,  \label{eq: eta + dynamics }
\end{align}
where 
\begin{align}
    v_1 &= \varepsilon\zs{{G\z{u^{L}, \mu^{L}} - G_{\textup{avg}}\z{u^{L}}}}, \label{eq: def v1}\\
    v_2 &= \varepsilon\zs{{G\z{u, \mu} - G_{\textup{avg}}\z{u} - G\z{u^{L}, \mu^{L}} + G_{\textup{avg}}\z{u^{L}}}}. \label{eq: def v2}
\end{align}

We use the superposition principle to determine the maximum value of $\eta$ by analysing the inputs $v_1$ and $v_2$ separately. Let us start with $v_1$. We append the subscripts to the notation of states to denote the time index. The discrete dynamics are given by:

\begin{align}
    \eta_{k + 1} &= (1 - \varepsilon)\eta_k + \varepsilon\left[G\z{u_k^{L}, \mu_k^{L}} - G_{\textup{avg}}\z{u_k^{L}}\right]. \label{eq: eta v1 def}
\end{align}
We define two additional variables:
\begin{align}
    \phi_{k + 1} &= \varepsilon\sum_{i = 1}^{k}\left[G\z{u_k^{L}, \mu_k^{L}} - G_{\textup{avg}}\z{u_k^{L}, \mu_k^{L}}\right] \label{eq: phi def} \\
    \theta_{k} &= \eta_{k} - \phi_{k} \label{eq: theta_def}
\end{align}

From \nref{eq: augmented system} and \nref{eq: theta_def} it holds
\begin{align}
    \theta_{k + 1} &= \eta_{k + 1} - \phi_{k + 1} = (1 - \varepsilon)\eta_k - \phi_k \red
    &= (1  - \varepsilon)(\eta_k - \phi_k) - (\varepsilon - 1)\phi_k - \phi_k \red
    &= (1 - \varepsilon) \theta_k - \varepsilon \phi_k \red
    &= (1 - \varepsilon)\zs{(1 - \varepsilon)\theta_{k - 1} - \varepsilon \phi_{k - 1}} - \varepsilon \phi_k \red
    &= (1 - \varepsilon)^2 \theta_{k - 1} - \varepsilon \zs{(1 - \varepsilon)\phi_{k - 1} + \phi_k}\red
    &= \dots \red
    &= -\varepsilon \zs{\sum_{i = 0}^{k}(1 - \varepsilon)^i \phi_{k - i}} \red
    &= \varepsilon \zs{\sum_{i = 0}^{k}(1 - \varepsilon)^i\z{\phi_k - \phi_{k - i}}}  - \varepsilon \zs{\sum_{i = 0}^k(1 - \varepsilon)^i \phi_k} \red
    &= \varepsilon \zs{\sum_{i = 0}^{k}(1 - \varepsilon)^i\z{\phi_k - \phi_{k - i}}}  - \phi_k \zs{1 - (1 - \varepsilon)^{k + 1}} \label{eq: theta def}
\end{align}
From \nref{eq: eta v1 def}, \nref{eq: phi def} and \nref{eq: theta def} we have
\begin{align}
    \eta_{k + 1}&= \theta_{k + 1} + \phi_{k + 1} \red
    &= {\varepsilon \zs{\sum_{i = 0}^{k}(1 - \varepsilon)^i\z{\phi_k - \phi_{k - i}}}} +  \z{\phi_{k + 1} - \phi_k} \red 
    &+ (1  -\varepsilon)^{k + 1}\phi_k. \label{eq: bounded inflation imperfect}
\end{align}
We use \nref{eq: concat average bound} in \nref{eq: bounded inflation imperfect} to derive:
\begin{align}
     \n{ \eta_{k + 1}} &\leq \underbrace{\varepsilon^2 \zs{\sum_{i = 0}^{k}(1 - \varepsilon)^i\z{i\, \sigma(L) + L\, \sigma(0)}}}_{S_1} + \varepsilon \sigma(L) \red 
     & + \varepsilon L\, \sigma(0) + \underbrace{\varepsilon (1 - \varepsilon)^{k + 1}\z{k\, \sigma(L) + L\, \sigma(0)}}_{S_2}
\end{align}
To compute $S_1$, we start by find the sum
\begin{align}
    \sum_{i = 1}^\infty i(1 - \varepsilon)^i &= \sum_{k = 1}^\infty \sum_{i = k}^\infty (1 - \varepsilon)^i = \sum_{k = 1}^\infty \frac{(1 - \varepsilon)^k}{\varepsilon} = \frac{1 - \varepsilon}{\varepsilon^2}.
\end{align}
Thus, we bound $S_1$ as follows:
\begin{align}
    S_1 \leq (1 - \varepsilon)\sigma(L) + \varepsilon L \sigma(0) \leq \sigma(L) + \varepsilon L \sigma(0).
\end{align}
For $S_2$, we define the function $z(x) = x\,(1 - \varepsilon)^x$. It is an easy exercises to check that the maximum of the function is given by $z\z{\frac{-1}{\log(1 - \varepsilon)}} = \frac{-e}{\log(1 - \varepsilon)}$. Therefore, for the bound of $S_2$ we have 
\begin{align}
    S_2 &\leq \varepsilon (1 - \varepsilon) \zs{\frac{-e}{\log \z{1 - \varepsilon}} \sigma(L) + L \sigma(0)}. \label{eq: bounded inflation semiperfect}
\end{align}
As $\lim_{\varepsilon \rightarrow 0^+} \tfrac{\varepsilon}{\log \z{1 - \varepsilon}} = 1$, for small enough $\varepsilon$, it follows:
\begin{align}
    &\leq e(1 - \varepsilon)\sigma(L) + \varepsilon(1 - \varepsilon) L\,\sigma(0) \red 
    &\leq e \sigma(L) + \varepsilon \, L \sigma(0).
\end{align}

Finally, we have
\begin{align}
     \n{\eta_{k + 1}} &\leq (1 + \varepsilon + e)\sigma(L) + 3\varepsilon\, L \sigma(0), \label{eq: bound v1}
\end{align}

which holds for all $k$. The norm can be made arbitrarily small by the right choice of parameters $L$ and $\varepsilon$. \\ \\
Now, we move on to the input $v_2$. We define the inflated boundary layer system:
\begin{align}
&\left.\begin{array}{l}
         {{u}_{\textup{bl}}^{\delta +}} \in u_{\textup{bl}}^\delta + \delta \mathbb{B}\\
         {{\mu}_{\textup{bl}}^{\delta+}} = M(u_{\textup{bl}}^\delta, \mu_{\textup{bl}}^\delta) 
    \end{array} \right\}, (u_{\textup{bl}}^\delta, \mu_{\textup{bl}}^\delta) \in \mathcal{U}  \times \Omega. \label{eq: boundary layer inflated}
\end{align} 
We claim the following:\\

\begin{lemma} \label{lemma: closeness of solutions}
For any period $L$, positive real number $a$ and compact set $K \in \mathcal{U}$, there exist a $\delta^*$ such that for every $\delta \in (0, \delta^*]$ and for any trajectory $(u_{\textup{bl}}^\delta, \mu_{\textup{bl}}^\delta)$ of the system in \nref{eq: boundary layer inflated} that is contained in $K \times \Omega$, there exist a concatenated trajectory $(u^L, \mu^L)$ such that 
\begin{align}
    \| G\z{u_{\textup{bl}}^\delta(k), \mu_{\textup{bl}}^\delta(k)} - G_{\textup{avg}}\z{u_{\textup{bl}}^\delta(k)} - \red G\z{u^{L}(k), \mu^{L}(k)} + G_{\textup{avg}}\z{u^{L}(k)}\| \leq a \kraj \label{eq: v2 bound}
\end{align}
\end{lemma}

\begin{pf}
$L, a, K$ are given. Let $\xi \coloneqq \col{u, \mu}$. Based on the continuity property of functions $G, G_{\textup{avg}}$, there exists $\rho > 0$ such that
\begin{align}
    \n{\xi_1 - \xi_2} \leq \rho \Rightarrow \n{G(\xi_1) - G(\xi_2)} \leq \tfrac{a}{2}, \red
    \n{G_{\textup{avg}}(\xi_1) - G_{\textup{avg}}(\xi_2)} \leq \tfrac{a}{2}. \label{eq: continuity boundness condition}
\end{align}
Next, we use \cite[Lemma 2]{wang2012analysis} for closeness of solutions of the inflated systems with parameters $(0, L + 1, \rho)$ and set $K \times \Omega$ to determine $\delta^*$. That means that for every trajectory of the system in \nref{eq: boundary layer inflated} where $\xi_{\textup{bl}}^\delta(k) \in K \times \Omega$ for all $k \in \dom(\xi_{\textup{bl}}^\delta)$, there exists a trajectory $\xi_{\textup{bl}}$ of the boundary layer system in \nref{eq: boundary layer system}, such that for each $k \in \dom(\xi_{\textup{bl}}^\delta)$ with $k \leq L + 1$, we have $\n{\xi_{\textup{bl}}^\delta(k) - \xi_{\textup{bl}}(k)} \leq \rho$. As the inflated boundary layer system is time invariant, any sample shifted trajectory is also a trajectory of the original system. Thus, for trajectories starting in $\xi_{\textup{bl}}^\delta(0), \xi_{\textup{bl}}^\delta(L), \dots, 
\xi_{\textup{bl}}^\delta(nL)$ with $n \in \N, nL \in \dom(\xi_{\textup{bl}}^\delta)$ there exist trajectories (not necessary the same one) $\xi_{\textup{bl}}$ such that the previous inequality holds for each segment of length $L$. We concatenate these boundary layer trajectories into $
\xi^L$ and write 
\begin{align}
    \n{\xi_{\textup{bl}}^\delta(k) - \xi^L(k)} \leq \rho \text{, for }k \in \dom(\xi_{\textup{bl}}^\delta). \label{eq: trajectory boundness condition}
\end{align}
From \nref{eq: continuity boundness condition} and \nref{eq: trajectory boundness condition} we conclude \nref{eq: v2 bound}. \end{pf}
The trajectories $\tilde{u}(k), {\mu}(k)$ of the transformed system in \nref{eq: transformed system}, where $\tilde{u}(k) \in K$ for all $k \in \dom\z{\tilde{u}}$ and $\n{\eta(0)} \leq a$, are also trajectories of the inflated boundary layer system in \nref{eq: boundary layer inflated} with 
\begin{align}
    \delta = \varepsilon\max_{u \in K, \mu \in \Omega, \eta \in K_{\eta}}\{G_{\textup{avg}}\z{u - \eta, \mu} - \eta\}, \label{eq: bound for delta}
\end{align}
 
where $K_{\eta}$ is the set in which $\eta$ is contained during the trajectory of the system. Let us prove that $K_{\eta} \subset a\mathbb{B}$ by first showing that i.e. $\n{\eta(1)} \leq a$. First we find $\varepsilon_1$ and $L$ such that $(1 + \varepsilon_1 + e)\sigma(L) + 3\varepsilon_1\, L \sigma(0) \leq \tfrac{a}{2}$. Then we use the same $L$, positive number $\tfrac{a}{2}$ and set $K$ with Lemma \ref{lemma: closeness of solutions} to find $\delta^*$. For $\varepsilon_2 = \frac{\delta^*}{\max_{u \in K, \eta \in a \mathbb{B}}\{G_{\textup{avg}}\z{u - \eta} - \eta\}}$, we guarantee that for one step, the solution of \nref{eq: transformed system} is also a solution of the inflated boundary layer system in \nref{eq: boundary layer inflated}. Thus, for $\varepsilon^* = \min\{\varepsilon_1, \varepsilon_2\}$, we have that variables in \nref{eq: def v1} and \nref{eq: def v2} are bounded as $\n{v_1} \leq \tfrac{a}{2}$ and $\n{v_2} \leq \tfrac{a}{2}$, and it follows from \nref{eq: eta + dynamics } that

\begin{align*}
    \n{\eta(1)} \leq \n{(1 - \varepsilon) \eta(0)} + \varepsilon a \leq   a.
\end{align*}
The next sample will also be a solution of the $\delta$-inflated boundary layer system and all of the previous bounds hold. Hence, the procedure can be repeated with the same $\delta^*$ for all $k \in \dom\z{\overline{u}}$, and it holds $\n{\varepsilon \eta(k)}\leq a$. \krajdokaz \end{pf}

Now, we return to the proof of Theorem \ref{thm: averaging theorem}. Let the set of initial conditions $K$ be given. From the stability of the set $\mathcal{A} \times \Omega$ in Assumption \ref{assum: stability of the average system} and the dynamics in \nref{eq: transformed system}, we have:

\begin{align}

    &V_{\textup{a}}(u^+ + \eta^+ + v^+, \mu^+) - V_{\textup{a}}(u + \eta + v, \mu)\red
    &=V_{\textup{a}}(\tilde{u}^+ + v^+, \tilde{\mu}) - V_{\textup{a}}(\tilde{u} + v, \tilde{\mu}) \red 
    &\leq V_{\textup{a}}(\tilde{u} + \varepsilon G_{\textup{avg}}(\tilde{u}, \tilde{\mu}) + v^+, \tilde{\mu}^+)  - V_{\textup{a}}(\tilde{u} + v, \mu) \red
    &- V_{\textup{a}}(\tilde{u} + \varepsilon G_{\textup{avg}}(\tilde{u}, \tilde{\mu}) + v^+, \tilde{\mu}^+) \red 
    &+ V_{\textup{a}}(\tilde{u} + \varepsilon G_{\textup{avg}}(\tilde{u} - \eta, \tilde{\mu}) - \varepsilon\eta + v^+)\red
    &\leq - \tilde{\alpha}_{\varepsilon}\z{{\varepsilon}}\alpha_{\textup{a}}\z{\n{\tilde{u} + v}_{\mathcal{A}}} + \varepsilon L_{V_{\textup{a}}}(1 + L_{\text{G}})\n{\eta}, \red
    &\quad\text{ for } \n{{\tilde{u} + v}}_{\mathcal{A}} \geq \alpha_{\gamma}\z{\gamma}\red
    &\leq -\hat{\alpha}_{\varepsilon}\z{{\varepsilon}}\alpha_{\textup{a}}\z{\n{\tilde{u} + v}_{\mathcal{A}}} \red
    &\quad\text{ for } \n{{\tilde{u} + v}}_{\mathcal{A}} \geq \max\{\alpha_{\gamma}\z{\gamma}, {\alpha}_{\varepsilon}\z{\varepsilon}\},
\end{align}
where $L_{\text{G}}$ and $L_{V_{\textup{a}}}$ are Lipschitz constants of the mapping $G_{\text{avg}}$ and function $V_{\text{a}}$ respectively, ${\alpha}_{\varepsilon}\z{\varepsilon} \geq \alpha_{\textup{a}}^{-1}\zs{ \frac{\varepsilon}{k\tilde{\alpha}_{\varepsilon}\z{{\varepsilon}}}L_{V_{\textup{a}}}(1 + L_{\text{G}})\n{\eta}}$, $\hat{\alpha}_{\varepsilon}\z{{\varepsilon}} \coloneqq(1 - k)\tilde{\alpha}_{\varepsilon}\z{{\varepsilon}}$ and $k \in (0, 1)$. The function $\alpha_\varepsilon$ is a function of class $\mathcal{K}$ on interval $(0, \bar{\varepsilon})$, as due to Assumption \ref{assum: stability of the average system}, $\tfrac{\varepsilon}{\tilde{\alpha}_\varepsilon(\varepsilon)}$ is bounded on that interval and $\eta$ can become arbitrarily small for proper choice of $\varepsilon$, per Lemma \ref{lemma: boundedness of eta}. Finally, we plug in the states of the original system to get 

\begin{align}
    &V_{\textup{a}}(u^+ + \eta^+ + v^+) - V_{\textup{a}}(u + \eta + v) \red 
    &\leq -\hat{\alpha}_{\varepsilon}\z{{\varepsilon}}\alpha_{\textup{a}}\z{\n{u + \eta + v}_{\mathcal{A}}} \red
    &\quad\text{ for } \n{{u + \eta + v}}_{\mathcal{A}} \geq \max\{\alpha_{\gamma}\z{\gamma}, {\alpha}_{\varepsilon}\z{\varepsilon}\}. \label{eq: proof average lyapunov diff end}
\end{align}
Let $\xi = u + \eta + v$. From the previous equation it follows 
\begin{align}
    &V_{\textup{a}}(\xi(k)) \leq V_{\textup{a}}(\xi(0)) -\sum_{i = 0}^{k - 1}\hat{\alpha}_{\varepsilon}\z{{\varepsilon}}\alpha_{\textup{a}}\z{\n{\xi(i)}_{\mathcal{A}}}\red
    &\quad\text{ for } \n{{\xi(k)}}_{\mathcal{A}} \geq \max\{\alpha_{\gamma}\z{\gamma}, {\alpha}_{\varepsilon}\z{\varepsilon}\}. \label{eq: proof lyapunov difference k}
\end{align}

Now, we move onto proving semi-global practical stability, Let $\Delta > \delta$ be any strictly positive real numbers. We choose parameters $\varepsilon$ and $\gamma$ such that $\oln{\eta}(\varepsilon) + \oln{v}(\gamma) \leq \tfrac{\delta}{4}$ and $\max\{\alpha_{\gamma}\z{\gamma}, {\alpha}_{\varepsilon}\z{\varepsilon}\} \leq \tfrac{\delta}{4}$. The conditional inequality in \nref{eq: proof lyapunov difference k} is satisfied when $\n{{u(k)}}_{\mathcal{A}} \geq \tfrac{\delta}{2}$. \\

\emph{\textbf{Semi-global stability}}\\ 

For ease of notation, we drop the explicit dependence on $\varepsilon$ and $\gamma$ in $\oln{\eta}(\varepsilon)$ and $\oln{v}(\gamma)$. We have to show that for any $R > \delta$, there exists $r > 0$, so that $\n{u(0)}_\mathcal{A} \leq r$ implies that $\n{u(k)}_\mathcal{A} \leq R$ for all $k \in \dom\z{u}$. From \nref{eq: assum lyapunov bound} and \nref{eq: proof lyapunov difference k}, it follows that
\begin{align}
    &\underline{\alpha}(\n{\xi(k)}_\mathcal{A}) \leq V_\textup{a}(\xi(k)) \leq V_\textup{a}(\xi(0)) \leq \oln{\alpha}(\n{\xi(0)}_\mathcal{A}) \red
    &\n{\xi(k)}_\mathcal{A} \leq \underline{\alpha}^{-1}\z{\oln{\alpha}(\n{\xi(0)}_\mathcal{A})} \red
    &\n{u(k)}_\mathcal{A} -\oln{\eta} - \overline{v} \leq \underline{\alpha}^{-1}\z{\oln{\alpha}(\n{u(0)}_\mathcal{A} + \oln{\eta} + \overline{v})} \red
    &\n{u(k)}_\mathcal{A} \leq \underline{\alpha}^{-1}\z{\oln{\alpha}(\n{u(0)}_\mathcal{A} + \oln{\eta} + \overline{v})} + \oln{\eta} + \overline{v} \red
    &\quad\text{ for } \n{{u}}_{\mathcal{A}} \geq \max\{\alpha_{\gamma}\z{\gamma}, {\alpha}_{\varepsilon}\z{\varepsilon}\} + \oln{\eta} + \oln{v}\nonumber
\end{align}

From last equation it follows that $R = \underline{\alpha}^{-1}\z{\oln{\alpha}(r + \oln{\eta} + \overline{v})} + \oln{\eta} + \overline{v}$. Thus, it holds $r = \oln{\alpha}^{-1}\z{\underline{\alpha}(R -\oln{\eta} - \overline{v})} -\oln{\eta} - \overline{v}$. Considering that the infimum value of $R$ is $\delta$ and that $\oln{\eta} + \overline{v} \leq \tfrac{\delta}{4}$, to ensure $r$ is positive, we assume $\varepsilon$ and $\gamma$ are chosen so that $\oln{\eta} + \overline{v} \leq \tfrac{1}{2} \oln{\alpha}^{-1}\z{\underline{\alpha}\z{\tfrac{3}{4}\delta}}$ which implies that $r \leq \tfrac{1}{2} \oln{\alpha}^{-1}\z{\underline{\alpha}\z{\tfrac{3}{4}\delta}}$. Furthermore, do assure that the Lyapunov difference is defined for those radiuses, we impose an additional inequality on the tuning parameters: $\max\{\alpha_{\gamma}\z{\gamma}, {\alpha}_{\varepsilon}\z{\varepsilon}\} + \oln{\eta} + \oln{v} <  \tfrac{1}{2} \oln{\alpha}^{-1}\z{\underline{\alpha}\z{\tfrac{3}{4}\delta}}$.\\ \\

\emph{\textbf{Practical attractivity}}\\ 
We have to show that for any $R, r$ that satisfy $\Delta > R > r > \delta > 0$, there exists $T$, such that $\n{u(0)}_\mathcal{A} \leq R$ implies that $\n{u(k)}_\mathcal{A} \leq r$ for all $k \in \dom\z{u}$ and $k \geq T$. First, we use the bound we derived in the proof of stability to define $r'\coloneqq \oln{\alpha}^{-1}\z{\underline{\alpha}(r -\oln{\eta} - \overline{v})} -\oln{\eta} - \overline{v}$, from which we can conclude that $\n{u(0)}_\mathcal{A} \leq r'$ implies that $\n{u(k)}_\mathcal{A} \leq r$ for all $k \in \dom\z{u}$. Let us define 
\begin{align}
    T \coloneqq \left\lceil\frac{\oln{\alpha}\z{R + \oln{\eta} + \overline{v}} - \underline\alpha\z{r' -\oln{\eta} - \overline{v}}}{\hat{\alpha}_\varepsilon\z{\varepsilon} \alpha \z{r -\oln{\eta} - \overline{v}}}\right\rceil + 1.
\end{align}
To prove via contradiction, we assume that $\n{u(k)} > r'$ for all $k \leq T$. Now, by using the upper and lower bound of the Lyapunov function on Equation in \nref{eq: proof lyapunov difference k}, it follows

\begin{align}
    &\underline{\alpha}\z{\n{u(k)}_\mathcal{A} -\oln{\eta} - \overline{v}} \leq \oln{\alpha}\z{R + \oln{\eta} + \overline{v}} \red
    &\quad \quad- k\hat{\alpha}_\varepsilon\z{\varepsilon}\alpha(r' -\oln{\eta} - \overline{v}) \red
    &\n{u(k)}_\mathcal{A} \leq  \underline\alpha^{-1}\z{\oln{\alpha}\z{R + \oln{\eta} + \overline{v}} - k\hat{\alpha}_\varepsilon\z{\varepsilon}\alpha(r' -\oln{\eta} - \overline{v})} \red
    &\quad \quad + \oln{\eta} + \overline{v} \label{eq: proof attractivity last equation}
\end{align}
Let us choose $k = T - 1$. When we plug in the chosen value of $k$ into inequality \nref{eq: proof attractivity last equation}, it follows that:
\begin{align}
    \n{u(k)}_\mathcal{A} \leq r',
\end{align}
which leads us to a contradiction. Thus, in the first $T$ steps, $u(k)$ trajectory will enter at least once the set $A + r'\mathbb{B}$. From the stability properties, we know that once the trajectory enters aforementioned set, it will never leave the set $\mathcal{A} + r \mathbb{B}$, which proves practical attractivity. \\ \\

Hence, to have semi-global practical asymptotic stability we have to choose our parameters $\varepsilon, \gamma$ so that they satisfy inequalities
\begin{align}
    &\oln{\eta}(\varepsilon) + \oln{v}(\gamma) \leq \tfrac{\delta}{4}\\
    &\max\{\alpha_{\gamma}\z{\gamma}, {\alpha}_{\varepsilon}\z{\varepsilon}\} \leq \tfrac{\delta}{4}\\
    &\max\{\alpha_{\gamma}\z{\gamma}, {\alpha}_{\varepsilon}\z{\varepsilon}\} + \oln{\eta} + \oln{v} <  \tfrac{1}{2} \oln{\alpha}^{-1}\z{\underline{\alpha}\z{\tfrac{3}{4}\delta}}
\end{align}
That concludes the proof of and semi-global practical asymptotic stability.\krajdokaz

\section{Proof of Theorem \ref{theorem: example1 convergence}}\label{proof: thm convergence full info}
First we show how to derive the boundary layer and averaged systems. Then we show that we can apply Theorem \ref{thm: averaging theorem} to prove stability. \\
The parameter $\alpha$ can be used as a time-scale separation parameter of the first layer in algorithm \nref{eq: example1 algorithm}. We derive the first boundary layer system ($\alpha = 0$):
\begin{align}
&\left\{ \begin{array}{ll}
     {\bfs{x}_{\textup{bl}}^1}^+ &= \bfs{x}_{\textup{bl}}^1 \\
     {\bfs{\xi}_{\textup{bl}}^1}^+ &= \bfs{\xi}_{\textup{bl}}^1 \\
     {\bfs{\mu}_{\textup{bl}}^1}^+ &= \mathcal{R}\bfs{\mu}_{\textup{bl}}^1
\end{array}\right., \label{eq: example1 bl system1}
\end{align}
and the first averaged system
\begin{align}
&\left\{ \begin{array}{ll}
    \hat{\bfs{x}}^+ &= (1 - \alpha\beta)\hat{\bfs{x}} + \alpha\beta \proj_{C}\z{\hat{\bfs{x}} - \gamma \bfs{\xi}} \\
    \hat{\bfs{\xi}}^+ &= (1 - \alpha)\hat{\bfs{\xi}} + \alpha (F(\hat{\bfs{x}}) + \mathcal{O}(\overline{a}))\\
    \hat{\bfs{\mu}}^+ &= \mathcal{R}\hat{\bfs{\mu}}
        \end{array}\right.,
\end{align}
which is an $\mathcal{O}(\oln{a})$ inflation of the nominal averaged system
\begin{align}
&\left\{ \begin{array}{ll}
        \hat{\bfs{x}}^+ &= (1 - \alpha\beta)\hat{\bfs{x}} + \alpha\beta \proj_{C}\z{\hat{\bfs{x}} - \gamma \hat{\bfs{\xi}}} \\
        \hat{\bfs{\xi}}^+ &= (1 - \alpha)\hat{\bfs{\xi}} + \alpha F(\hat{x})\\
        \hat{\bfs{\mu}}^+ &= \mathcal{R}\hat{\bfs{\mu}}
        \end{array}\right. .\label{eq: example1 nominal average system}
\end{align}
Furthermore, we use $\alpha \beta$ for the parameter of the second-time layer separation to determine the second boundary layer system
\begin{align}
&\left\{ \begin{array}{ll}     
        {\bfs{x}_{\textup{bl}}^2}^+ &= \bfs{x}_{\textup{bl}}^2\\
        {\bfs{\xi}_{\textup{bl}}^2}^+ &= (1 - \alpha)\bfs{\xi}_{\textup{bl}}^2 + \alpha F(\bfs{x}_{\textup{bl}}^2)\\
        {\bfs{\mu}_{\textup{bl}}^2}^+ &= \mathcal{R}{\bfs{\mu}_{\textup{bl}}^2}
     \end{array}\right.,\label{eq: example1 bl system2}
\end{align}
and the second averaged system
\begin{align}

&\left\{ \begin{array}{ll}
     \tilde{\bfs{x}}^+ &= (1 - \alpha\beta)\tilde{\bfs{x}} + \alpha\beta \proj_{C}\z{\tilde{\bfs{x}} - \gamma F(\tilde{\bfs{x}})}\\
     \tilde{\bfs{\xi}}^+ &= (1 - \alpha)\tilde{\bfs{\xi}} + \alpha F(\tilde{\bfs{x}})\\
     \tilde{\bfs{\mu}}^+ &= \mathcal{R}\tilde{\bfs{\mu}}
         \end{array}\right.,\label{eq: example1 average system2}
\end{align}
which is the algorithm in \nref{eq: example1 fb algorithm} with additional bounded dynamics that renders the set $\{\bfs{x}^*\} \times C_{\textup{F}} \times \mathbb{S}^m$ UGAS. \\
In order to satisfy Assumption \ref{assum: boundness of the average aprox} for both averaged systems, we establish the following result:
\begin{lemma} \label{lemma: sinusoidal aproximation}
    For any solution of the first boundary layer system $(\bfs{x}_{\textup{bl}}^1, \bfs{\xi}_{\textup{bl}}^1, \bfs{\mu}_{\textup{bl}}^1)$ and compact set $C$ such that $\bfs{x}_{\textup{bl}}^1 \in C$ for all $k \in \dom\z{\bfs{x}_{\textup{bl}}^1} $, it holds that:
\begin{align}
    &\Biggl \| \frac{1}{N}\sum_{i = 1}^{N}\Bigl [ {2}{A^{-1}} J\z{\bfs{x}_{\textup{bl}}^1(i) + A\mathbb{D}\bfs{\mu}_{\textup{bl}}^1(i)}\mathbb{D}\bfs{\mu}_{\textup{bl}}^1(i) \red
    &- F(\bfs{x}_{\textup{bl}}^1(i)) - \mathcal{O}(\bar{a})\Bigr ] \Biggr \| \leq \sigma_1\z{N}, \label{eq: example average bound system 1}
\end{align} 
where $\sigma_1: \R^+ \rightarrow \R^+$ is a function of class $\mathcal{L}$.\kraj
\end{lemma}
\begin{pf}
See Appendix \ref{proof: lemma example1 frist average}. \krajdokaz
\end{pf}
\begin{lemma}\label{lemma: example1 first average}
    For any solution of the second boundary layer system $(\bfs{x}_{\textup{bl}}^2, \bfs{\xi}_{\textup{bl}}^2)$ and compact set $C$ such that $\col{\bfs{x}_{\textup{bl}}^1, \bfs{\xi}_{\textup{bl}}^1}  \in C$ for all $k \in \dom\z{\bfs{x}_{\textup{bl}}^1} $, it holds that:
\begin{align}
    \n{\frac{\gamma}{N}\sum_{i = 1}^{N}\zs{\bfs{\xi}_{\textup{bl}}^2(i) - F(\bfs{x}_{\textup{bl}}^2(i))}} \leq \sigma_2\z{N}, \label{eq: example1 average bound system 2}
\end{align} 
where $\sigma_2: \R^+ \rightarrow \R^+$ is a function of class $\mathcal{L}$.\kraj
\end{lemma}\label{lemma: example1 second average}
\vspace{-0.6cm}
\begin{pf}
See Appendix \ref{proof: lemma example1 second average}. \krajdokaz
\end{pf}
\vspace{-0.6cm}
To prove stability, we start from the second layer and ``move" upwards. As the second averaged system satisfies Assumptions \ref{assum: boundness of the average aprox} due to Lemma \ref{lemma: sinusoidal aproximation}, Assumption \ref{assum: continuity of jump and averag functions} due the nonexpansivnes of the projection mapping \cite[Prop. 12.28, 29.1]{bauschke2011convex}, Assumption \ref{assum: boundedness of hidden dynamics} due to Lemma \ref{lemma: boundedness of eta} and Assumption \ref{assum: stability of the average system} due to \nref{eq: example1 algorithm lyapunov}, we have that due to Theorem \ref{thm: averaging theorem}, the nominal averaged system in \nref{eq: example1 nominal average system} renders the set $\{x^*\} \times C_{\textup{F}} \times \mathbb{S}^m$ SGPAS as $\alpha\beta \rightarrow 0$, with the Lyapunov difference given by

\begin{align}
    &V(\tilde{\bfs{x}}^+ + \eta_1^+) - V(\tilde{\bfs{x}} + \eta_1) \leq \red
    &\quad - \tfrac{1}{2}\alpha\beta(1 - c)(2 - \alpha\beta c)\n{\tilde{\bfs{x}} + \eta_1 - \bfs{x}^*}_\mathcal{A}\red
    &\text{for }\n{\tilde{\bfs{x}} + \eta_1}_\mathcal{A} \geq \alpha_\varepsilon(\alpha\beta) \geq \sqrt{\tfrac{\alpha\beta\, L \n{\eta_1}}{2\alpha\beta(1 - c)(2 - \alpha\beta c)}}\label{eq: example1 second layer average system lyapunov}
\end{align}

and the perturbation dynamics 

\begin{align}
    \eta_1^+ = (1 - \alpha\beta)\eta_1 + \alpha\beta (\proj_{C}\z{\tilde{\bfs{x}} - \gamma F(\tilde{\bfs{x}})} - \tilde{\bfs{\xi}}). \label{eq: example 1 bounded dynamics eta second averge system}
\end{align}
We note that we had to take $\alpha\beta$ as the time-scale separation parameter. If we had chosen only $\beta$, as might be the intuition, the function $\alpha_\varepsilon$ of class $\mathcal{K}$ that appears in the inequality in \nref{eq: example1 second layer average system lyapunov}, would have an implicit dependence on the parameter $\alpha$. In fact, decreasing $\alpha$ would increase the value of the function $\alpha_\varepsilon$, as it would hold
\begin{align}
    \alpha_\varepsilon(\beta) \geq \sqrt{\tfrac{\beta\, L \n{\eta_1}}{2\alpha\beta(1 - c)(2 - \alpha\beta c)}},
\end{align}
which would invalidate all of the following stability analysis. Thus, it is important to capture all of the parameters that affect the speed of convergence. Nevertheless, if we assume that the parameter $\alpha$ is contained in the set $(0, \oln{\alpha})$, it is possible to construct a function of class $\mathcal{K}$, such that it holds\, 
$\n{\tilde{x} + \eta_1}_\mathcal{A} \geq \alpha_\beta(\beta) \geq \alpha_\varepsilon(\alpha\beta)$. Hence, the averaged system in \nref{eq: example1 nominal average system} renders the set $\{x^*\} \times C_{\textup{F}} \times \mathbb{S}^m$ SGPAS as $\beta \rightarrow 0$.\\
The first averaged system is an $\mathcal{O}(\overline{a})$ inflation of the nominal averaged system, and it can be shown that the inflation introduces a small perturbation in the Lyapunov difference inequality which can be made arbitrarily small by choosing $\oln{a}$ small enough. For the sake of the proof, we set $\oln{a} = \alpha_a(\beta)$, where $\alpha_a$ is a function of class $\mathcal{K}$. Thus, it also satisfies Assumption \ref{assum: stability of the average system}, Assumption \ref{assum: boundness of the average aprox} due to \nref{eq: example 1 bounded dynamics eta second averge system}, Assumption \ref{assum: continuity of jump and averag functions} because of \cite[Prop. 12.28, 29.1]{bauschke2011convex}, and Assumption \ref{assum: boundedness of hidden dynamics} as a result of \nref{eq: example 1 bounded dynamics eta second averge system}. Hence, the system in \nref{eq: example1 algorithm} renders the set $\{x^*\} \times C_{\textup{F}} \times \mathbb{S}^m$ SGPAS as $(\alpha, \oln{a}, \beta) \rightarrow 0$.

\krajdokaz

\section{Proof of Lemma \ref{lemma: sinusoidal aproximation}}\label{proof: lemma example1 frist average}

Without the loss of generality, let a solution of the first boundary layer system is given by $\bfs{x}_{\textup{bl}}^1(k) = \bfs{x}_{\textup{bl}}^1 = const.$, $\bfs{\xi}_{\textup{bl}}^1 = const.$, $\mu_{\textup{bl}}^{1,i}(k) = \col{(\sin(\omega_i^j k), \cos(\omega_i^j k))_{j \leq m_i}}$, $\mu_{\textup{bl}}^1(k) = \col{(\mu_{\textup{bl}}^{1,i}(k)))_{i \in \mathcal{I}}}$. First, with the following Lemma, we characterize the properties of average discrete-time sinusoidal signals: 
\begin{lemma}
For any $\phi, \phi_i, \phi_l \in \R$ such that $\phi \neq 2\pi t$,  $\phi_i \pm \phi_l \neq  2\pi p, t, p \in \Z$, it holds that:

\begin{align}
    &\frac{1}{N}\n{\sum_{k = 0}^{N - 1}\sin(\phi k)} \leq \frac{c_1}{N}, \label{eq: example1 sin single bound} \\
    &\frac{1}{N}\n{\sum_{k = 0}^{N - 1}\cos(\phi k)} \leq \frac{c_2}{N}, \label{eq: example1 cos single bound} \\
    &\frac{1}{N}\n{\sum_{k = 0}^{N - 1}\sin(\phi_i k) \sin(\phi_l k)} \leq \frac{c_3}{N}, \label{eq: example1 sin bound} \\
    &\frac{1}{N}\n{\sum_{k = 0}^{N - 1}(\sin^2(\phi k) - \tfrac{1}{2})} \leq \frac{c_4}{N}, \label{eq: example1 sin square bound}
\end{align}
for some $c_1, c_2, c_3, c_4 > 0$.\kraj
\end{lemma}
\begin{pf}
We note that for $\phi \neq 2\pi p, p\in \Z$, it follows:
\begin{align}
    &  \n{\sum_{k = 0}^{N-1} e^{j\phi k}}^2 =  \n{\sum_{k = 0}^{N-1} \cos(\phi k)}^2 + \n{\sum_{k = 0}^{N-1} \sin(\phi k)}^2 \red
    &= \n{\frac{e^{j\phi N} - 1}{e^{j\phi} - 1}}^2 \leq \frac{4}{\n{e^{j\phi} - 1}^2}.
\end{align}
Therefore, we have
\begin{align}
    \n{\sum_{k = 0}^{N-1} \cos(\phi k)} \leq \frac{2}{\n{e^{j\phi} - 1}} \coloneqq c_1, \label{eq: example1 cos only bound}\\
    \n{\sum_{k = 0}^{N-1} \sin(\phi k)} \leq \frac{2}{\n{e^{j\phi} - 1}} \coloneqq c_2.
\end{align}
Equations \nref{eq: example1 sin single bound}, \nref{eq: example1 cos single bound} follow from the previous equations.\\
Let $\phi = \phi_i \pm \phi_l$. From \nref{eq: example1 cos only bound}, we have

\begin{align}
    &\n{\sum_{k = 0}^{N-1}\z{\sin(\phi_ik)\sin(\phi_lk) \mp \cos(\phi_ik)\cos(\phi_lk) }} \red
    &\leq \frac{2}{\min_{\phi\in  \phi_i \pm \phi_l}\n{e^{j\phi} - 1}} \coloneqq c_3\label{eq: example1 sine sum}
\end{align}

For any scalars $a, b \in \R, c \in \R_+$, that satisfy equation $\n{a \pm b} \leq c$, it holds
\begin{align}
    -c \leq a + b &\leq c \text{, and} \red
    -c \leq a - b &\leq c. \label{eq: example1 absolute value condition}
\end{align}
By summing the last two inequalities, we have
\begin{align}
    \n{a} \leq c. \label{eq: example1 absolute value result}
\end{align}
Thus from \nref{eq: example1 sine sum}, \nref{eq: example1 absolute value condition} and \nref{eq: example1 absolute value result}, we conclude
\begin{align}
    \n{\sum_{k = 0}^{N-1}{\sin(\phi_ik)\sin(\phi_lk)}} \leq c_3.
\end{align}
Again, \nref{eq: example1 sin bound} follows trivially. Finally, using the identity $\sin^2(x) = \frac{1 - \cos(2x)}{2}$, we rewrite Equation \nref{eq: example1 sin square bound} as

\begin{align}
    \frac{1}{2N}\n{\sum_{k = 0}^{N - 1}\cos(2\phi k)}  \leq \frac{c_4}{N}.
\end{align}
By switching $2\phi$ instead of $\phi$ in \nref{eq: example1 cos only bound}, analogously it is possible to prove \nref{eq: example1 sin square bound}.\krajdokaz
\end{pf}

Via the Taylor expansion of the an addend in \nref{eq: example average bound system 1}, we have
\begin{align}
&{2}{A^{-1}} J\z{\bfs{x}_{\textup{bl}}^1 + A\mathbb{D}\mu_{\textup{bl}}^1(i)}\mathbb{D}\mu_{\textup{bl}}^1(i) -F(\bfs{x}_{\textup{bl}}^1) - \mathcal{O}(\bar{a}) \red
&= {{2}A^{-1}} J(\bfs{x}_{\textup{bl}}^1)\mathbb{D}\mu_{\textup{bl}}^1(i) + F(\bfs{x}_{\textup{bl}}^1)^\top \mathbb{D}\mu_{\textup{bl}}^1(i) \mathbb{D}\mu_{\textup{bl}}^1(i) \red 
&- F(\bfs{x}_{\textup{bl}}^1)
\end{align}

Due to the inequality $\sqrt{\sum_{i = 0}^m x_i^2} \leq \sum_{i = 0}^m|x_i|$, we can bound the expression in \nref{eq: example average bound system 1} via the bounds for each row:
\begin{align}
     &\Biggl \| \frac{1}{N}\sum_{i = 1}^{N}\Bigl [ {2}{A^{-1}} J\z{\bfs{x}_{\textup{bl}}^1 + A\mathbb{D}\mu_{\textup{bl}}^1(i)}\mathbb{D}\mu_{\textup{bl}}^1(i) \red
    &- F(\bfs{x}_{\textup{bl}}^1) - \mathcal{O}(\bar{a})\Bigr ] \Biggr \| \leq \red
    &\sum_{j = 1}^m \left\| \sum_{k = 0}^{N-1} \left[ \tfrac{2}{a_j}J_j(\bfs{x}_{\textup{bl}}^1) \sin(\phi_j k) \right. \right. \red
    & + \nabla_{x_j}J_j(\bfs{x}_{\textup{bl}}^1)(2\sin^2(\phi_j k) - 1) \red
    & + 2\sum_{l \neq j}^n \nabla_{x_l}J_j(\bfs{x}_{\textup{bl}}^1) \sin(\phi_l k)\sin(\phi_j k) \Bigr ] \Biggr \| \leq \red
    &\frac{2m}{\underline{a}}\n{J(\bfs{x}_{\textup{bl}}^1)} \frac{c_1}{N} +  \n{\nabla J(\bfs{x}_{\textup{bl}}^1)}_{\infty}\frac{2\overline{c}\,m^2}{N},
\end{align}
where $\underline{a} \coloneqq \min_i a_i$ and $\overline{c} \coloneqq \max \{c_3, c_4\}$. Thus for the compact set $C$, we define $\sigma_1(N) \coloneqq \frac{2m}{\underline{a}}\max_{x \in C}\n{J(x)} \frac{c_1}{N} +  \max_{x \in C}\n{\nabla J(x)}_{\infty}\frac{2\overline{c}\,m^2}{N}$ which belongs to the class of $\mathcal{L}$ functions. \krajdokaz

\section{Proof of Lemma \ref{lemma: example1 first average}}\label{proof: lemma example1 second average}
A solutions of the second boundary layer system are given by $\bfs{x}_{\textup{bl}}^2(k) = \bfs{x}_{\textup{bl}}^2 =  const.$, and $\bfs{\xi}_{\textup{bl}}^2(k) = (1 - \alpha)^k \z{\bfs{\xi}_{\textup{bl}}^2(0) - F(\bfs{x}_{\textup{bl}}^2)} + F(\bfs{x}_{\textup{bl}}^2)$. Thus, the norm in \nref{eq: example1 average bound system 2} can be rewritten as 
\begin{align}
     \n{\frac{\gamma}{N}\sum_{i = 0}^{N - 1}\zs{\bfs{\xi}_{\textup{bl}}^2(i) - F(\bfs{x}_{\textup{bl}}^2(i))}} \leq \red \frac{\gamma}{N}\n{\sum_{i = 0}^{N - 1}\zs{(1 - \alpha)^i\z{\bfs{\xi}_{\textup{bl}}^2(0) - F(\bfs{x}_{\textup{bl}}^2})}} \leq \red
     \frac{\gamma}{N\alpha} \n{{{\bfs{\xi}_{\textup{bl}}^2(0) - F(\bfs{x}_{\textup{bl}}^2)}}}.
\end{align}
Thus, for the compact set $C$, we define $\sigma_2(N) \coloneqq \frac{\gamma}{N\alpha} \max_{\col{\bfs{\xi}, \bfs{x}} \in C}\n{{{\bfs{\xi} - F(\bfs{x})}}}$, which belongs to the class of $\mathcal{L}$ functions. \krajdokaz
 
 \section{Proof of Theorem \ref{thm: essentisally cyclic forward algorithm}}\label{appendix: proof essentisally cyclic forward algorithm}
For notational simplicity, we denote $S(k) \coloneqq \hat{S}_x\z{k; \bfs{\tau}_0}$. Thus, the algorithm reads as
\begin{align}
    &\left\{ \begin{array}{ll}
    \bfs{x}^+&= \bfs{x} - \alpha {S}\z{k}  F(\bfs{x}) \\
    k^+ &= k + 1.
    \end{array}\right. \label{eq: proof example 2 reduced discrete time alg}
\end{align}

One epoch is defined as $r$ iterations of the algorithm in \nref{eq: proof example 2 reduced discrete time alg}, where $r$ is the period of the function $S$ from Lemma \ref{lemma: example 2 periodicity}. From the proof of the Lemma, it follows that every agent \emph{individually} jumps $r_i$ times in one epoch. Let $F_j: \R^m \rightarrow \R^m$ be the mapping that returns the rows of the pseudogradient that correspond to the agents that sample at $k = j + rn$, $n \in \N$, $0 \leq j \leq r - 1$, i.e. $F_j(\bfs{x}) \coloneqq S(j)F(\bfs{x})$. We define the full update operator, the asynchronous update operator, and error operator respectively, as

\begin{align}
    &T \coloneqq I - \alpha \Gamma F, \\
    &E \coloneqq (I - \alpha F_1)\dots(I - \alpha F_r), \\
    &R \coloneqq \tfrac{1}{\alpha}(T - E).
\end{align}
where $\Gamma \coloneqq \diag{(r_i I_{m_i})_{i\in \mathcal{I}}}$. We note that the $E$ operator represents one epoch of the algorithm in \nref{eq: proof example 2 reduced discrete time alg}, i.e. $\bfs{x}(k + r) = E(x(k))$. \\
The proof of convergence is analogous to the proof in \cite{chow2017cyclic}, and here we just provide the outlines. The error operator can be bounded as 

\begin{align}
\n{R(x)}^2 \leq \frac{\alpha^2 L^2 r^2 \oln{r}^2}{2}(1 + \alpha L)^{2m}\n{S(x)}^2,
\end{align}
where $\oln{r} = \max_{i \in \mathcal{I}} r_i$. For the Lyapunov function candidate, we propose
\begin{align}
    V(\bfs{x}) = \n{\bfs{x} - \bfs{x}^*}_{\Gamma^{-1}}^2.
\end{align}
It can be proven that
\begin{align}\label{proof: epoch iteration}
\n{E(\bfs{x}) - \bfs{x}^*}_{\Gamma^{-1}}^2 \leq \z{1  -\frac{\alpha \mu_{\textup{F}}^2\underline{r}^2}{2}}\n{\bfs{x} - \bfs{x}^*}_{\Gamma^{-1}}^2, 
\end{align}
where $\underline{r} = \min_{i \in \mathcal{I}} r_i$, if $\alpha$ is chosen such that

\begin{align}
    \tfrac{1}{2} + \tfrac{\oln{r}\eta}{\mu_\textup{F}^2} - (1 - \alpha\eta)\z{1 - \alpha - \tfrac{\z{\alpha L r \oln{r} (1 + \alpha L)^r }^2}{2\eta}} \leq 0. \nonumber
\end{align}
The inequality is satisfied for $\eta$ and $\alpha$ small enough, and $ \alpha \ll \eta$. We note that the inequality does not depend on the initial conditions  of the timers $\bfs{\tau}_0$. Equation \nref{proof: epoch iteration} holds for \emph{epochs}, not necessarily the individual samples. Due to Lipschitz continuity of the pseudogradient, it follows that 

\begin{align}
    \n{\bfs{x}^+ - \bfs{x}^*}_{\Gamma^{-1}} \leq (1 + \alpha L\tfrac{\oln{r}}{\underline{r}})\n{\bfs{x} - \bfs{x}^*}_{\Gamma^{-1}}.\nonumber
\end{align}
Thus, for some $k = ir + j$, where $i, j \in \N$ we have
\begin{align}
    &\n{\bfs{x}(k) - \bfs{x}^*}_{\Gamma^{-1}}^2 \leq \red
    &\z{1 + \alpha L\tfrac{\oln{r}}{\underline{r}}}^{2r}\z{1 -\frac{\alpha \mu_{\textup{F}}^2\underline{r}^2}{2}}^i\n{\bfs{x}(0) - \bfs{x}^*}_{\Gamma^{-1}}^2.
\end{align}
It holds $i = \lfloor\tfrac{k}{r}\rfloor \geq \tfrac{k}{r} - 1$. Hence, the previous inequality becomes

\begin{align}
    &\n{\bfs{x}(k) - \bfs{x}^*}_{\Gamma^{-1}}^2 \leq \red
    &\z{1 + \alpha L\tfrac{\oln{r}}{\underline{r}}}^{2r}\z{1 -\frac{\alpha \mu_{\textup{F}}^2\underline{r}^2}{2}}^{\tfrac{k}{r} - 1}\n{\bfs{x}(0) - \bfs{x}^*}_{\Gamma^{-1}}^2.\nonumber
\end{align}
The last inequality is the KL exponential stability bound for all initial conditions of timer states $\bfs{\tau}_0$. Thus, the dynamics in \nref{eq: proof example 2 reduced discrete time alg} render $\bfs{x}^*$ UGES. \\ \\
Additionally, we need to establish the Lyapunov difference convergence speed. To do this, we construct a Lyapunov function using a similar procedure as in the proof of \cite[Thm. 4.14]{khalil2002nonlinear}, which we omit due to space limitations. Let 
\begin{align}
    a &\coloneqq \frac{\z{1 + \alpha L\tfrac{\oln{r}}{\underline{r}}}^{2r}}{1 -\frac{\alpha \mu_{\textup{F}}^2\underline{r}^2}{2}} \red
    b & \coloneqq \z{1 -\frac{\alpha \mu_{\textup{F}}^2\underline{r}^2}{2}}^{\tfrac{1}{r}}.
\end{align}
Then, the Lyapunov function satisfies the following properties

\begin{align}
    &\frac{1 - \z{1 - \alpha L\tfrac{\oln{r}}{\underline{r}}}^{2\delta\,r}}{2\alpha L\tfrac{\oln{r}}{\underline{r}}(1 - \alpha L\tfrac{\oln{r}}{\underline{r}})}\n{\bfs{x} - \bfs{x}^*}_{\Gamma^{-1}}^2 \leq V(\bfs{x}) \red 
    &\leq \frac{a(b^\delta\,r - 1)}{b - 1}\n{\bfs{x} - \bfs{x}^*}_{\Gamma^{-1}}^2 \red 
    &V(\bfs{x}^+) - V(\bfs{x}) \leq - (1 - ab^{\delta\,r})\n{\bfs{x} - \bfs{x}^*}_{\Gamma^{-1}}^2, \nonumber
\end{align}  
where $\delta$ is large integer. Using the Taylor expansion, for small values of $\alpha$, it holds

\begin{align}
    (1 - ab^{\delta\,r}) \approx \alpha \z{\frac{\mu^2 \underline{r}^2 (\delta - 1)}{2} - 2rL\tfrac{\oln{r}}{\underline{r}}}.
\end{align}
Thus, for $\delta$ large enough, we can guarantee that the Lyapunov difference is negative. Furthermore, we have that $\tfrac{\alpha}{(1 - ab^{\delta\,r})}$ is bounded on some interval $(0, \alpha^*)$ and the Lyapunov function satisfies the conditions from Assumption \nref{assum: stability of the average system}.\krajdokaz 
\vspace{-0.4cm}
\section{Proof of Lemma \ref{lemma: example 2 periodicity}} \label{appendix: proof of (T,q) period}
\vspace{-0.4cm}

First, let us denote least common sampling time as $T = T_i \tfrac{p}{p_i}$. \\

\emph{Claim 1.} The number of jumps in any time interval $[\hat{t}, \hat{t} + T)$ is constant. \kraj \\
Let us denote as the number of jumps in this interval as $q$. We ``slide" the interval by some $\Delta T$, i.e. $[\hat{t} + \Delta T, \hat{t} + T + \Delta T)$, so that we exclude one event in $[\hat{t}, \hat{t} + \Delta T)$. As $T = T_i r_i$ for all $i\in\mathcal{I}$, it follows that there must be an event in the interval $[\hat{t} + T, \hat{t} + T + \Delta T)$, thus the total number of jumps in the interval $[\hat{t} + \Delta T, \hat{t} + T + \Delta T)$ remains the same. We can repeat this procedure for any $\Delta \hat{T}$ by sliding the interval for every jump by $\Delta T_j$ until $\Delta \hat{T} = \sum_j \Delta T_j$. The number of jumps in the interval is equal to $q$.\\ 

\emph{Claim 2.} For every $(\hat{t}, \hat{j}) \in \dom(\bfs{\tau})$ where a jump occurred, it holds $\bfs{\tau}(\hat{t}, \hat{j}) = \bfs{\tau}(\hat{t} + T, \hat{j} + r)$. \kraj\\
We observe that if a jump is initiated by an agent at $t = \hat{t}$, these same agent will also initiate a jump at $t = \hat{t} + T$. Furthermore, if an agents did not jump at $\hat{t}$, it will also not jump at $\hat{t} + T$.  Thus, in the moment $t = \hat{t} + T$, the same agent will jump, i.e. $\tau_i = T_i$ and it follows that $\bfs{\tau}(\hat{t}, \hat{j}) = \bfs{\tau}(\hat{t} + T, \hat{j} + r)$.\\

As the functions $S_{x}$ and $S_{\tau}$ are single-valued,   the claim of the Lemma holds.\krajdokaz
\vspace{-0.4cm}
\section{Proof of Theorem \ref{theorem: example2 convergence}}
\vspace{-0.4cm}
As the proof is analogous to the proof of Theorem \ref{theorem: example1 convergence}, we provide just the required system definitions and averaging Lemmas. The equivalent discrete-time system of \nref{eq: example2 algorithm} is given by

\begin{align} 
    &\left\{ \begin{array}{ll}
        \bfs{x}^{+} &= \bfs{x} - \alpha \beta S(k)\bfs{\xi}\\
        \bfs{\xi}^{+} &= \bfs{\xi} + \alpha S(k)\z{2A^{-1} J(\bfs{x} + A\mathbb{D}\bfs{\mu})\mathbb{D}\bfs{\mu} - \bfs{\xi}} \\
        \bfs{\mu}^{+} &= (I - S_{\mu}(k))\bfs{\mu} + S_{\mu}(k))\mathcal{R}\bfs{\mu} \\ 
        \bfs{\kappa}^+ &= \bfs{\kappa} + S_\tau(\bfs{\tau}) \\
        k^+ &= k + 1.
            \end{array}\right. 
\end{align} 
The first boundary-layer system is defined as
\begin{align} 
    &\left\{ \begin{array}{ll}
        {\bfs{x}^1_{\textup{bl}}}^{+} &= \bfs{x}^1_{\textup{bl}} \\ 
         {\bfs{\xi}^1_{\textup{bl}}}^{+} &= \bfs{\xi}^1_{\textup{bl}} \\
         {\bfs{\mu}^1_{\textup{bl}}}^{+} &= (I - S_{\mu}(k))\bfs{\mu}^1_{\textup{bl}} + S_{\mu}(k))\mathcal{R}\bfs{\mu}^1_{\textup{bl}} \\
         {\bfs{\kappa}}^+ &= \bfs{\kappa} + S_\tau(k) \\
         k^+ &= k + 1,
             \end{array}\right. \label{eq: example 2 first bl system}
\end{align}

while the first averaged system is given as
\begin{align} 
    &\left\{ \begin{array}{ll}
        \tilde{\bfs{x}}^{+} &= \tilde{\bfs{x}} - \alpha \beta S(k)\tilde{\bfs{\xi}}\\ 
        \tilde{\bfs{\xi}}^{+} &= \tilde{\bfs{\xi}} + \alpha S(k)\z{F(\tilde{\bfs{x}}) - \tilde{\bfs{\xi}}} \\ 
        \tilde{\bfs{\mu}}^{+} &= (I - S_{\mu}(k))\tilde{\bfs{\mu}} + S_{\mu}(k))\mathcal{R}\tilde{\bfs{\mu}} \\
        {\bfs{\kappa}}^+ &= \bfs{\kappa} + S_\tau(k) \\
        k^+ &= k + 1.
         \end{array}\right. 
\end{align} 
The second boundary-layer system follows the dynamics
\begin{align} 
    &\left\{ \begin{array}{ll}
        {\bfs{x}^2_{\textup{bl}}}^{+} &= \bfs{x}^2_{\textup{bl}} \\
         {\bfs{\xi}^2_{\textup{bl}}}^{+} &= \bfs{\xi}^2_{\textup{bl}} - \alpha S(k)(F(\bfs{x}^2_{\textup{bl}}) - \bfs{\xi}^2_{\textup{bl}})\\
         {\bfs{\mu}^2_{\textup{bl}}}^{+} &= (I - S_{\mu}(k))\bfs{\mu}^2_{\textup{bl}} + S_{\mu}(k))\mathcal{R}\bfs{\mu}^2_{\textup{bl}} \\
         {\bfs{\kappa}}^+ &= \bfs{\kappa} + S_\tau(k)\\
         k^+ &= k + 1,
         \end{array}\right.
\end{align} 
whereas the second averaged system is defined as
\begin{align} 
&\left\{ \begin{array}{ll}
        \tilde{\bfs{x}}^{+} &= \tilde{\bfs{x}} - \alpha \beta S(k)F(\tilde{\bfs{x}})\\
        \tilde{\bfs{\xi}}^{+} &= \tilde{\bfs{\xi}} + \alpha S(k)\z{F(\tilde{\bfs{x}}) - \tilde{\bfs{\xi}}} \\
        \tilde{\bfs{\mu}}^{+} &= (I - S_{\mu}(k))\tilde{\bfs{\mu}} + S_{\mu}(k))\mathcal{R}\tilde{\bfs{\mu}} \\
        {\bfs{\kappa}}^+ &= \bfs{\kappa} + S_\tau(k) \\
        k^+ &= k + 1.
         \end{array}\right.
\end{align} 

To prove Assumption \ref{assum: boundness of the average aprox}, the following two Lemmas are needed:

\begin{lemma}\label{lemma: example 2 first average}
    For any solution of the first boundary layer system $(\bfs{x}_{\textup{bl}}^1, \bfs{\xi}_{\textup{bl}}^1, \bfs{\mu}_{\textup{bl}}^1)$ and compact set $C$, such that $\bfs{x}_{\textup{bl}}^1 \in C$ for all $k \in \dom\z{\bfs{x}_{\textup{bl}}^1} $, it holds that:
\begin{align}
    &\Biggl \| \frac{1}{N}\sum_{i = 1}^{N}\Bigl [ {2}{A^{-1}} J\z{\bfs{x}_{\textup{bl}}^1(i) + A\mathbb{D}\bfs{\mu}_{\textup{bl}}^1(i)}\mathbb{D}\bfs{\mu}_{\textup{bl}}^1(i) \red
    &- F(\bfs{x}_{\textup{bl}}^1(i)) - \mathcal{O}(\bar{a})\Bigr ] S(i) \Biggr \| \leq \sigma_1\z{N}, \label{eq: lemma example 2 boundness of approx for sinus}
\end{align} 
where $\sigma_1: \R^+ \rightarrow \R^+$ is a function of class $\mathcal{L}$.\kraj
\end{lemma}
\begin{pf}
First, we note that the difference in the previous inequality for rows corresponding to agent $i$ is equal to zero whenever the agent is not jumping. This motivates us to study the ``isolated" system of agent $i$, instead of the group dynamics in \nref{eq: example 2 first bl system}. Consider the dynamics

\begin{align}
    \mu^i(\kappa_i + 1) = \mathcal{R}_i \mu^i(\kappa_i).
\end{align}
Without the loss of generality, let a solution of the previous system be given by $\mu^i(\kappa_i) = \operatorname{col}((\cos(\omega_i^j \kappa_i), $ $\sin(\omega_i^j \kappa_i))_{j\leq m_i})$. The solution of $\mu_{\textup{bl}}^{1,i}$ is similar to the solution of $\mu^i$, as it also has the same samples, but they ``persist" for more iterations, i.e. until the agent $i$ jumps again. If we define the set valued mapping $k(v, i) \coloneqq \{u \mid \kappa_i(u) = v\}$ that relates the global jump counter to the internal counter of agent $i$, it holds that $\mu_{\textup{bl}}^{1,i}(k(v, i)) = \mu^i(v)$. Furthermore, at the sample $v$ given by the internal counter for agent $i$, or given by $\oln{k}(v, i) = \max k(v, i)$ by the global counter, the internal counter of some other agent $j$ is given by

\begin{align}\label{proof: representation of other samples}
    \kappa_j(v, i) \coloneqq \left\lfloor \Delta_j^i + \tfrac{T_i}{T_j}v\right\rfloor.
\end{align}
where $\Delta_j^i\coloneqq  \kappa_j^0 + \tfrac{\tau_i - \kappa_i^0 T_i}{T_j}$. Hence it holds $\mu_{\textup{bl}}^{1,j}(\oln{k}(v, i)) = \mu^j(\kappa_j(v, i))$. Lastly, we see that diagonal elements of $S(k)$ corresponding to agent $i$ are different from zero for $k \in \{u \mid u = \bar{k}(v, i), v \in \N\}$. Let us denote the norm in \nref{eq: lemma example 2 boundness of approx for sinus} as $\Phi$ and write $N = r\, l + o, l \in \N, 0\leq o < r$. The previous iterator relations and properties of $S(k)$ allows us bound the inequality as

\begin{align}
    &N\Phi \leq \sum_{i = 1}^M\Biggl \|\sum_{v = \kappa_i(0)}^{\kappa_i(0) + r_i\,l - 1}\Bigl [ - F(\bfs{x}_{\textup{bl}}) - \mathcal{O}(\bar{a}) + \red
    &\tfrac{2}{a_i} J_i\z{\bfs{x}_{\textup{bl}}^1 + A\mathbb{D}\bfs{\mu}_{\textup{bl}}^1(\oln{k}(v, i)}\mathbb{D}_i{\mu}_{\textup{bl}}^{1, i}(\oln{k}(v, i))  \Bigr ] \Biggr \| + \Phi_\textup{r} \red
        
    &\leq \sum_{i = 1}^M \Biggl \| \sum_{v = \kappa_i(0)}^{\kappa_i(0) + r_i\,l - 1}\Bigl [ \tfrac{2}{a_i} J_i\z{\bfs{x}_{\textup{bl}}^1}\mathbb{D}_i \mu_i(v)  \red 
    & + 2 \nabla_{x_i}J_i(\bfs{x}_{\textup{bl}}^1)^\top \mathbb{D}_i \mu_i(v) \mathbb{D}_i \mu_i(v) -  \nabla_{x_i}J_i(\bfs{x}_{\textup{bl}}^1)  \red 
    &+ 2\sum_{j \neq i}^M \nabla_{x_j}J_i(\bfs{x}_{\textup{bl}}^1)^\top \mathbb{D}_j \mu_i(\kappa_j(v, i)) \mathbb{D}_i \mu_i(v)\Bigr ] \Biggr \| + \Phi_\textup{r} \red
    &\leq \sum_{i = 1}^M \sum_{u = 1}^{m_i} \Biggl \| \sum_{v = \kappa_i(0)}^{\kappa_i(0) + r_i\,l - 1}\Bigl [\tfrac{2}{a_i} J_i\z{\bfs{x}_{\textup{bl}}^1}\sin\z{\omega_i^u v} \red
    & \nabla_{x_i^u}J_i(\bfs{x}_{\textup{bl}}^1)\zs{2 \sin\z{\omega_i^u v} \sum_{j = 1}^{m_i}\sin(\omega_i^j v)-1} \red
    & + 2\sum_{j \neq i}^M\sum_{s = 1}^{m_j} \nabla_{x_j^s}J_i(\bfs{x}_{\textup{bl}}^1)\sin(\omega_j^s\kappa_j(v, i)) \sum_{j = 1}^{m_i}\sin(\omega_i^j v)\Bigr ] \Biggr \| \red & + \Phi_\textup{r},     
\end{align}
where

    \begin{align}
    &\Phi_\textup{r} \coloneqq \sum_{i = 1}^M\,r_i\, m\, M\, \max_{\bfs{x} \in C}\n{\tfrac{1}{a_i}J\z{\bfs{x} + A\mathbb{B}}} + \n{F(\bfs{x}) + \mathcal{O}(\oln{a})}.
\end{align} 

Using Lemma \nref{lemma: example1 first average} and Assumption \ref{assum: example 2 different freq}, we can derive the upper bounds of all of the sums in the norm, apart for the last one, which contains addends of  form $\sin(\omega_i^j v)\sin(\omega_j^s\kappa_j(v, i))$. Using the same procedure as in proof of Lemma \ref{lemma: sinusoidal aproximation} and Equation \nref{proof: representation of other samples}, we find the equivalent exponential representation. For some $\omega_1, \omega_2$, consider the sum 

\begin{align}
    &\sum_{v = \kappa_i(0)}^{\kappa_i(0) + r_i\,l - 1}e^{j\z{\omega_1v + \omega_2\left\lfloor \Delta_j^i + \tfrac{p_i}{p_j}v\right\rfloor}} \red
    &=\sum_{v = \kappa_i(0)}^{\kappa_i(0) + r_i - 1} e^{j\z{\omega_1 v + \omega_2\left\lfloor \Delta_j^i + \tfrac{p_i}{p_j}v\right\rfloor}}\times \sum_{u = 0}^{l - 1}e^{j\z{\omega_1r_i + \omega_2 r_j}u} \red
    &=\frac{1 - e^{j\z{\omega_1r_i + \omega_2 r_j}l}}{1 - e^{j\z{\omega_1r_i + \omega_2 r_j}}}\sum_{v = \kappa_i(0)}^{\kappa_i(0) + r_i - 1} e^{j\z{\omega_1 v + \omega_2\left\lfloor \Delta_j^i + \tfrac{p_i}{p_j}v\right\rfloor}}, \nonumber
\end{align}
where the second equality follows from Assumption \ref{assumption example2 sampling ratio}, the properties of the least common multiple $p = p_i\, r_i = p_j\, r_j$, and last equality holds for $\omega_1 r_i + \omega_2 r_j \neq 2\pi z, z\in \Z$. Thus, we have 
\begin{align}
    \n{\sum_{v = \kappa_i(0)}^{\kappa_i(0) + r_i\,l - 1}e^{j\z{\omega_1v + \omega_2\left\lfloor \Delta_j^i + \tfrac{p_i}{p_j}v\right\rfloor}} }^2 \leq \red
    \n{\frac{r_i}{1 - e^{j\z{\omega_1r_i + \omega_2 r_j}}}}^2\leq  C_e, \nonumber
\end{align}
where $C_e$ is the supremum with respect to all possible combinations of $r_i$ and $\omega_i^j$. The rest of the procedure follows the same steps as after Equation \nref{eq: example1 sine sum} in the proof of Lemma \ref{lemma: example1 first average}. The bound that holds regardless of the initial conditions of the timers $\bfs{\tau}(0)$. Thus, the Lemma holds. \krajdokaz
\end{pf}
\begin{lemma}\label{lemma: example 1 second average}
For any solution of the second boundary layer system $(\bfs{x}_{\textup{bl}}^2, \bfs{\xi}_{\textup{bl}}^2)$, it holds that:

\begin{align}
     \n{\frac{\gamma}{N}\sum_{i = 1}^{N}S(i)\zs{\bfs{\xi}_{\textup{bl}}^2(i) - F(\bfs{x}_{\textup{bl}}^2(i))}} \leq \sigma_2\z{N},
\end{align}

where $\sigma_2: \R^+ \rightarrow \R^+$ is a function of class $\mathcal{L}$.\kraj
\end{lemma}
\begin{pf}
The proof is analogous to the proof of Lemmas \ref{lemma: example1 first average} and \ref{lemma: example 2 first average} and is omitted due to space constraints. \krajdokaz
\end{pf}

\end{document}